\documentclass[12pt]{article}

\usepackage{amsfonts}
\usepackage{amsmath}
\usepackage{amssymb}
\usepackage{fullpage}
\usepackage{bbold}
\usepackage{graphicx}
\usepackage{float}
\usepackage[caption=false]{subfig}
\usepackage[font=small,labelfont=bf]{caption}

\newcommand{\EllipticK}{\mathbf{K}}
\newcommand{\EllipticE}{\mathbf{E}}

\title{Dynamical integrity of the safe basins in a problem of forced escape}

\author{Pavel Kravetc\thanks{Faculty of Mechanical Engineering, Technion -- Israel Institute of Technology, Haifa, Israel}\and Oleg Gendelman$^*$\and Alexander Fidlin\thanks{Institute of Engineering Mechanics, Karlsruhe Institute of Technology, Karlsruhe, Germany}}

\begin{document}
\maketitle

\abstract{This paper explores the use of the Approximation of Isolated Resonance (AIR) method for determining the safe basins (SBs) in the problem of escape from a potential well. The study introduces a novel approach to capture the location and the shape of the SBs and establish their erosion profiles. The research highlights the concept of "true" safe basins, which remain invariant with phase shifts, a critical factor often faced in real-world applications. A cubic polynomial potential serves as the benchmark to illustrate the proposed method.}

\section{Introduction}

The concept of escape from a potential well is an essential topic that finds recurrent appearances across various scientific and engineering fields. Its extensive range of applications encompasses numerous phenomena, such as the capsizing of ships~\cite{Virgin1989, belenky2007}, energy harvesting~\cite{mann2009energy}, the pull-in phenomenon in micro-electro-mechanical systems (MEMS)~\cite{5482087, leus2008dynamic}, chemical reactions~\cite{kramers1940, fleming1993activated}, the buckling of arches~\cite{VIRGIN1992357, champneys2019happy}, and the physics behind the Josephson junction~\cite{barone1982physics}, among others.

To study the escape of a particle captured by a primary resonance due to the effect of external forcing Approximation of Isolated Resonance (AIR) method was proposed in~\cite{Gendelman2018}. The AIR method typically proceeds in the following manner. Initially, the equations of motion are converted into a more manageable format using action-angle (AA) variables. Assuming the presence of a primary $1:1$ resonance, a slow phase can be chosen, allowing for averaging over all the rapid phases. This process yields a slow-flow equation that depicts the resonance manifold (RM). Assuming the absence of damping, this slow-flow movement can be represented by a non-trivial conservation law. By analyzing the resulted slow-flow system, one can establish the critical forcing amplitude necessary for the escape as well as the safe basins (SBs), i.e., all the initial conditions that correspond to non-escaping trajectories.

Approximation of isolated resonance was first applied to safe basins in~\cite{Karmi2021}.  The exploration of safe basins for a weakly-damped particle within a truncated quadratic potential well (i.e., level-crossing problem) was discussed in~\cite{attila}. Escape in finite time and the corresponding SBs were studied in~\cite{10.1063/5.0142761}.

In order to quantify the size of the SB the concept of the \emph{integrity measure} was introduced~\cite{doi:10.1098/rspa.1989.0009}. The primary types of the integrity measures include the Global Integrity Measure (GIM)~\cite{doi:10.1098/rspa.1989.0009}, the Local Integrity Measure (LIM)~\cite{soliman1989integrity} and the Integrity Factor (IF)~\cite{REGA2005902}.

The Global Integrity Measure (GIM) represents the hyper-volume of the safe basin. In simpler terms, it is the total volume (or area in 2D problems) of the space that represents safe initial conditions. However, the GIM might not always be a prudent measure as it includes all parts of the safe basin, including the fractal ones~\cite{doi:10.1098/rspa.1989.0009, soliman1989integrity}.

On the other hand, if the SB is a basin of attraction (BOA), one can define the Local Integrity Measure (LIM) as the normalized minimum distance from the attractor (the point or set towards which the system evolves) to the boundary of the safe basin. Note that, LIM is a property of the attractor and while it lacks clear theoretical background it is often a convenient choice for numerical assessment of dynamical integrity of the BOA~\cite{Habib2021}. A more general approach suitable for potential wells is the Integrity Factor (IF) defined as the radius of the largest hyper-sphere entirely belonging to the safe basin. Both LIM and IF measures focus only on the compact core of the safe basin, excluding the fractal parts~\cite{soliman1989integrity, lenci2003optimal}.

These integrity measures are used to study the loss of integrity when system parameters change. The process of the system losing its integrity is sometimes referred to as erosion, and the changes in integrity across different system parameters form what is known as the erosion profile~\cite{REGA2005902}. Analyzing erosion profiles is crucial in safe engineering design~\cite{lenci2019global, rega2021global}.

In real life applications it is imperative to get results that are valid for a range or parameters, as often we do not know the exact values of the parameters. For example, in many instances the phase of the external excitation remains undetermined. Hence, the notion of true safe basins was introduced in~\cite{regaLenci2008}. The ``true" SB is the safe basin which is invariant to the phase shift. It can be defined as an intersection of the safe basins for every value of the phase.

In this work the authors continue exploring the application of AIR method to study the safe basins. In particular, we propose a method to capture the location and the shape of the SBs and establish their erosion profiles. We choose a cubic polynomial as a benchmark potential to illustrate the method.

The paper is organized as follows. Section~\ref{sec:main} introduces all the necessary notation and describes the AIR method to establish the safe basins including the ``true" safe basins. Section~\ref{sec:cubic} contains application of the AIR method to an a cubic potential. Finally, Section~\ref{sec:conclusion} presents conclusions and discussions.

\section{Main Idea}\label{sec:main}
Consider a single-degree-of-freedom classical particle of the unit mass trapped inside a local minimum of some potential~$V(q)$. Assume, the particle is subject to the influence of external harmonic forcing. Then, the equation of motion is
 \begin{equation}\label{eq:main}
 \ddot{q} + \frac{d V}{d q} = F \sin{\left(\Omega t + \psi\right)},
 \end{equation}
 where~$q$ denotes the displacement, dot represents derivative with respect to the time~$t$, and~$F$,~$\Omega$ and~$\psi$ are the amplitude, the frequency and the phase of the external forcing, respectively. We assume the damping to be negligible. In fact, while taking damping into account is crucial in many engineering applications, a small damping does not drastically affect the analysis~\cite{10.1063/5.0142761}, therefore, we limit the scope of the present work to the undamped system.

We define the escape from potential well when the total energy~$E=p^2/2 + V(q)$ of the particle exceeds some predefined threshold value~$E_\text{thres}$, i.e.,
\[
\max\limits_t\left\{E(t)\right\}\ge E_\text{thres}.
\] Here, $p = \dot{q}$ denotes the momentum of the particle.

The safe basin is defined as a union of all initial conditions that correspond to non-escaping trajectories. If~$D_\psi$ denote a SB for a particular value~$\psi$ of the external forcing phase, then the ``true" safe basin~$D_{\text{true}}$ is defined as an intersection of $D_\psi$ for all values of $\psi \in \left[0,\;2\pi\right)$. 

In~\cite{Karmi2021} the authors attempted to locate and quantify the safe basins using isolated resonance approximation (AIR) method. Recall, there are two mechanisms of escape based on the structure of the resonance manifold: the maximum mechanism and the saddle mechanism. The escape through the maximum mechanism occurs when the trajectory on the RM passes tangentially to the line~$\xi=\xi_{\max}$. The reason for that is the limiting trajectory passing in the proximity of the separatrix of the basic hamiltonian $H_o(q, p) = E$. Isolated resonance approximation disregards the dynamics at the separatrix tangle, hence, the fails to deliver accurate predictions. On the other hand, the method allows to utilize a truncated potential, and if an appropriate truncation level~$\xi_\text{max}$ is used, it's possible to identify the safe basin while eliminating the unsafe fractal tongues. Note that, escape from a truncated potential is equivalent to a level-crossing problem~\cite{attila}. In order to find $\xi_\text{max}$, one can utilize estimates of the width of the separatrix chaotic layer.

A typical framework of the AIR method goes as follows. First, it is necessary to obtain a slow-flow equation: averaged motion of the particle on the resonance manifold. In order to do so, one needs to rewrite the equation using action-angle (AA) variables and perform averaging over all the fast phases (under the assumption of the primary resonance). Second, one can study the dynamics on the resonance manifold and identify the safe basins. Finally, by applying the transformation back to $(q, p)$ variables, one can approximate the safe basins on the phase plane.

Note that, in the absence of damping, i.e., when $\lambda=0$, the slow-flow motion admits a non-trivial first integral. Assume the following equation describes the RM:
\begin{equation}
C(\vartheta, J) = C_o, \qquad \text{for}\quad C_o \in \mathbb{R},
\end{equation} where 
\begin{equation}\label{eq:slow:phase:definition}
\vartheta = \theta - \omega t -\psi
\end{equation} is a slow phase and~$J$ denotes averaged action. The constant term~$C_o$ is defined by the initial condition at which the system is captured by the resonance. Then, one can identify a safe basin~$D$ with its boundary~$\partial D$ defined as a phase curve on the RM for a particular value~$C_o^*$ of the constant~$C_o$:
\begin{equation}
\partial D = \left\{C(\vartheta, J) = C_o^*\right\}.
\end{equation}

Using the change of variables $\Psi: (\vartheta, J) \mapsto (q, p)$, one can find the safe basin~$\Psi(D)$ on the $(q,p)$-plane. Note that, the transformation~$\Psi$ is canonical, hence, area preserving, therefore
\[
\text{area}(\Psi(D)) = \text{area}(D).
\]

In many cases depending on the form of the potential in question, it is impossible to obtain the integral of motion with the averaged action~$J$. In this case, one can express slow motion using averaged energy~$\xi$ instead. Transformation~$\Psi$ is not canonical anymore, therefore the area of the defined safe basin~$\Psi(D)$ is
\begin{equation}
\text{area}(\Psi(D))=\int\limits_{\Psi(D)}\text{d}q\,\text{d}p = \int\limits_{D}L(\vartheta, \xi)\;\text{d}\vartheta\,\text{d}\xi,
\end{equation} where
\begin{equation}
L(\vartheta, \xi) = \det{\left(\text{J}\Psi\right)} = \frac{\partial q}{\partial \vartheta}\frac{\partial p}{\partial \xi} - \frac{\partial p}{\partial \vartheta}\frac{\partial q}{\partial \xi},
\end{equation} where $\text{J}$ denotes the Jacobian.

Recall, in~\cite{Karmi2021} the SBs are classified based on their geometry on the $(\vartheta, \xi)$-cylinder. In particular, we distinguish safe basins of maximum type (SBMT) and safe basins of saddle type (SBST). SBMT's boundary is defined as a level set that passes tangentially at the escaping threshold~$\xi=\xi_{\max}$, i.e., it is defined as
\begin{equation}\label{eq:sbmt}
C(\vartheta,\; \xi) = C(\vartheta_{\max},\; \xi_{\max}),
\end{equation} where~$\vartheta_{\max}$ is the value of~$\vartheta$ variable at the tangent point. SBMTs can be further classified into SBMT of the first kind and SBMT of the second kind based on its topology: SBMT of the second kind ``wraps" around the entire $(\vartheta, \xi)$-cylinder, while SBMT of the first kind does not. Note that SBMT of kind I is bounded by~\eqref{eq:sbmt} and the bottom of the RM cylinder.

The boundary of the SBST is defined as
\[
C(\vartheta,\; \xi) = C(\vartheta_{\text{saddle}},\; \xi_{\text{saddle}}), \qquad \xi < \xi_{\text{saddle}} 
\] where $\left(\vartheta_{\text{saddle}},\; \xi_{\text{saddle}}\right)$ is the saddle point.

Keep in mind that two variants of SB can exist simultaneously. However, this dual presence is only feasible for SBMT II and SBST. When both SBMT I and SBST are present on the $(\vartheta, \xi)$-cylinder, SBST is entirely included within SBMT I. Consequently, the splitting of SB occurs during the transition from SBMT I to SBMT II, which precisely takes place at the saddle connection:
\begin{equation}\label{eq:saddle:connection}
C(\vartheta_{\max},\; \xi_{\max}) = C(\vartheta_{\text{saddle}},\; \xi_{\text{saddle}}). 
\end{equation} By solving equation~\eqref{eq:saddle:connection} for~$F$ one can find the critical value $\widehat{F}$ when the safe basin splits into two components. Note that the split of the SB for a particular value~$\psi$ into two components corresponds to the disappearance of the ``true" SB.

Note that by the definition~\eqref{eq:slow:phase:definition}, the phase~$\psi$ enters the slow phase~$\vartheta$ linearly, hence it represents a horizontal rotation of the RM on the $(\vartheta, \xi)$-cylinder. If $D_\psi$ denotes a safe basin for a particular value~$\psi$, one can find the ``true" safe basin~$D_{\text{true}}$ as an intersection of $D_\psi$ for all the value of $\psi\in\left[0,2\pi\right)$. Therefore, it is sufficient to find the minimum value~$\widetilde{\xi}$ on the boundary~$\partial D_\psi$ and then the ``true" safe basin boundary~$\partial D_{\text{true}}$ is just a circle~$\xi=\widetilde{\xi}$. Then, the cylinder $[0,2\pi)\times[0, \widetilde{\xi}]$ defines the ``true" safe basin. The area of the ``true" SB in $(q,p)$-plane is simply expressed as:
\[
\text{area}(D) = 2\pi J(\widetilde{\xi}),
\] where~$J$ denotes the action.

%

\section{Cubic Potential}\label{sec:cubic}
Consider the system~\eqref{eq:main} with a cubic potential
\begin{equation}\label{eq:cubic}
V(q) = \frac{q^2}{2}-\frac{q^3}{3}.
\end{equation}
Isolated resonance approximation for the system with the potential~\eqref{eq:cubic} was used in~\cite{FARID2020105182} to establish critical forcing amplitude~$F_{\text{crit}}(\Omega)$ as a function of the excitation frequency~$\Omega$. In this section, we apply the same framework to locate and quantify safe basins for a particular value of the phase~$\psi$ as well as the ``true" SB.

Using a slow phase~$\vartheta$ and averaged energy~$\xi$ as variables one can write the slow motion on the $1:1$ resonance manifold as the following conservation law:
\begin{equation}\label{eq:cons:cub}
C(\vartheta, \xi) = \xi - \frac{F \pi^2 \sqrt{3} \sin{(z)}}{k^2 \EllipticK^2(k)} \, \frac{Q}{1-Q^2} \cos{\vartheta} - \Omega J(\xi) = \text{const},
\end{equation} where 
\begin{gather*}
z = \frac{1}{3}\arccos{(1-12\xi)},\qquad Q = \exp{\left(-\frac{\pi \EllipticK(k^\prime)}{\EllipticK(k)}\right)},
\end{gather*}
and the averaged action $J(\xi)$ is
\[
J(\xi) = \frac{2\sqrt{2}}{3^{1/4} 5 \pi} \sqrt{\sin{\left(\frac{2\pi}{3}-z\right)}}\left(\frac{3}{2} \EllipticE(k) - \sqrt{3} \sin{\left(\frac{\pi}{3}-z\right)}\cos{(z)}\EllipticK(k)\right).
\]
Functions~$\EllipticK(k)$ and~$\EllipticE(k)$ denote the complete elliptic integrals of the first and the second kind, respectively, with the modulus
\[
k = \sqrt{\frac{\sin{z}}{\sin{(\frac{2\pi}{3}-z)}}},
\] and $k^\prime$ denotes the complementary modulus, i.e., $k^\prime=\sqrt{1-k^2}$.

The displacement~$q$ and the momentum~$p$ are expressed as functions of angle~$\vartheta$ and energy~$\xi$ as follows:
\begin{align}
q(\vartheta, \xi) &= q_{\text{min}} + \left(q_\text{max} - q_\text{min}\right)\text{sn}^2\left(\frac{\vartheta \EllipticK(k)}{\pi}, \, k\right),\label{eq:q:trans}\\
p(\vartheta, \xi) &= \sqrt{\frac{2}{3}} \sqrt{c-q_{\min}} \;(q_{\max}-q_{\min})\cdot\notag\\ &
   \text{cn}\left(\frac{\theta  \EllipticK(k)}{\pi }, \, k\right)
   \text{dn}\left(\frac{\theta  \EllipticK(k)}{\pi }, \, k\right)
   \text{sn}\left(\frac{\theta  \EllipticK(k)}{\pi }, \, k\right),\label{eq:p:trans}
\end{align} where $\text{cn}$, $\text{sn}$ and $\text{dn}$ are the Jacobi elliptic functions. For the detailed derivation of the formulae~\eqref{eq:cons:cub}--\eqref{eq:p:trans} we refer the reader to~\cite{FARID2020105182}.

Finally, functions $q_{\min}(\xi) < q_{\max}(\xi) < c(\xi)$ denote the three roots of the cubic equation~$V(q) = \xi$ for $0<\xi<1/6$. Using trigonometric formulae, one can express them as follows:
\begin{align*}
q_{\min}(\xi)&=\frac{1}{2}-\sin \left(\frac{1}{3} \arccos (1-12 \xi )+\frac{\pi }{6}\right),\\
q_{\max}(\xi)&=\frac{1}{2}-\sin \left(\frac{\pi }{6}-\frac{1}{3} \arccos(1-12 \xi )\right), \\
c(\xi)&=\cos \left(\frac{1}{3} \arccos(1-12 \xi )\right)+\frac{1}{2}.
\end{align*}

Consider a curve~$C(\vartheta, \xi) = C_{\text{saddle}}$, where $C_{\text{saddle}} = C(\vartheta=0, \xi=\xi_{\text{saddle}})$ and $\xi_{\text{saddle}}$ is the value of~$\xi$ at the saddle point, i.e., $\xi_{\text{saddle}}$ is a solution to the following equation:
\[
\frac{\partial}{\partial \xi} C(0, \xi) = 0.
\]

One can define a function~$\xi_\text{s}(\vartheta)$ by solving equation~$C(\vartheta, \xi) = C_{\text{saddle}}$ for~$\xi < \xi_\text{saddle}$. This curve defines the boundary of the SBST on the RM for any value of the phase~$\psi$. The ``true" safe basin can be obtained by taking a circle $\xi = \widehat{\xi}$ passing tangentially through the point of the curve~$C(\vartheta, \xi) = C_{\text{saddle}}$ where the value of $\xi$ achieves the minimum, i.e., $\widehat{\xi}$ is the solution to $C(\pi, \xi) = C_{\text{saddle}}$. In other words, it is~$\widehat{\xi} = \min\limits_\vartheta{\xi_\text{s}(\vartheta)} = \xi_{\text{s}}(\pi)$.

Similarly, one can obtain the boundary curve of SBMT:
\[
C_{\max} = C(\pi, \xi_{\max})= C_{\max}.
\] Recall, if SBMT is of the II kind the values of angle~$\vartheta$ span the entire period from 0 to $2\pi$. Therefore, the~$\xi$ component of the boundary curve can be expressed as a function~$\xi_m(\vartheta)$ with the minimum value achieved at~$\vartheta = 0$. Therefore,~$\widehat{\xi} = \min_{\vartheta}\{\xi_m\} = \xi_m (0)$.

Recall, that the ``true" safe basin of the maximum type disappears abruptly when parameter~$F$ crosses a critical point~$\widehat{F}$. It happens precisely when the boundary curve becomes the saddle connection. Therefore, we can find the exact value of~$\widehat{F}$ by solving the following equation:
\begin{equation}
C(\pi, \xi_{\text{max}})=C(0, \xi_{\text{saddle}}).
\end{equation}

\begin{figure}[H]
\centering
\subfloat[$F=0.015$]{\includegraphics*[width=0.3\textwidth]{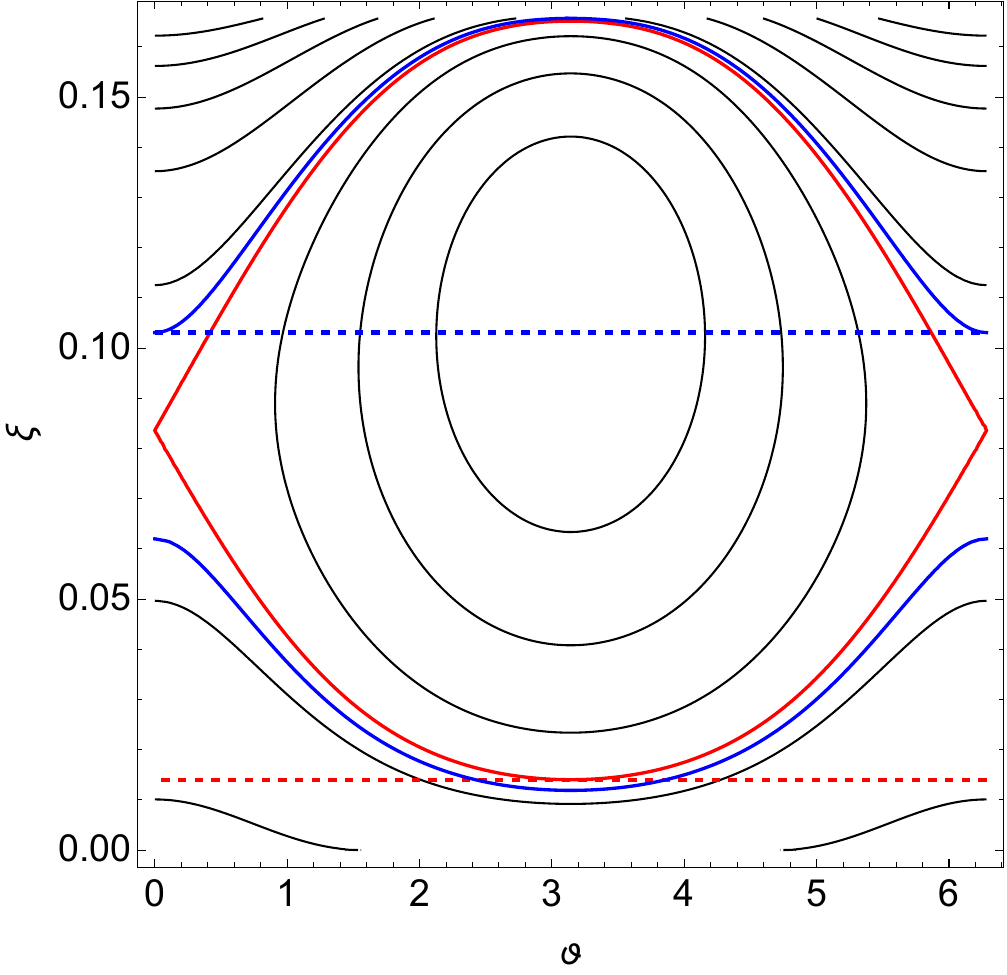}}
\hspace{10pt}
\subfloat[$F=\widehat{F}\approx 0.0155721$]{\includegraphics*[width=0.3\textwidth]{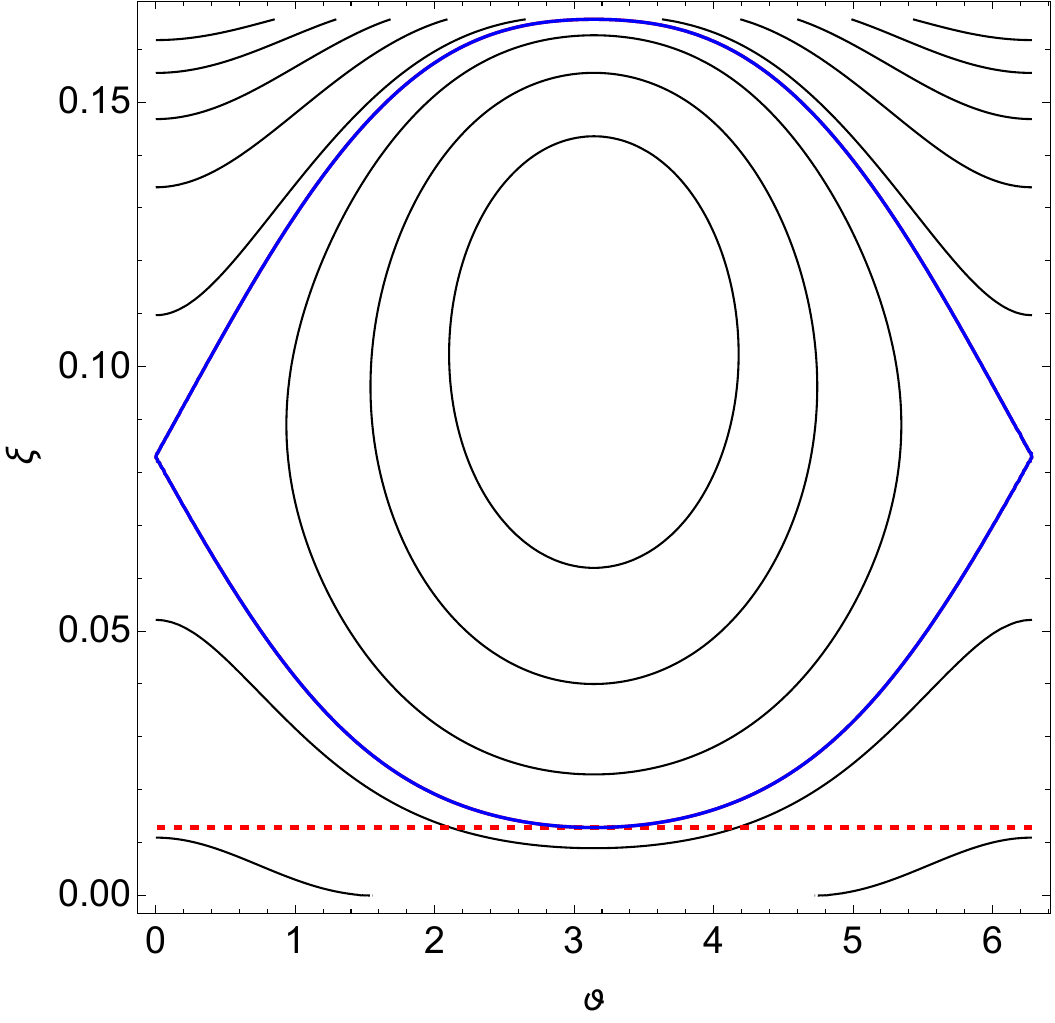}}
\hspace{10pt}
\subfloat[$F=0.016$]{\includegraphics*[width=0.3\textwidth]{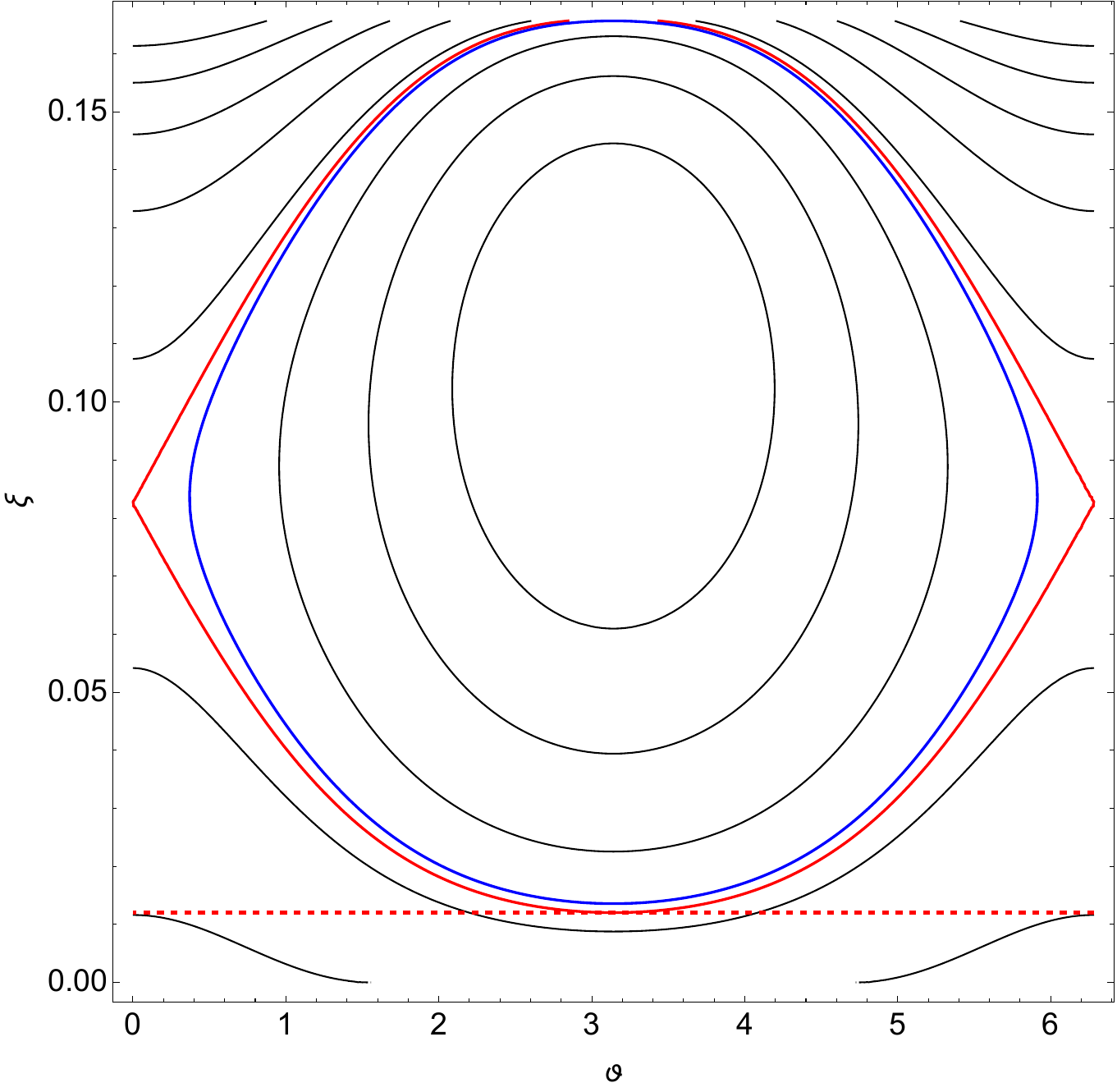}}
\caption{Passing through the saddle-connection on the $(\vartheta,\xi)$-cylinder. Blue and red curves are the level sets defined by the maximum and the saddle mechanisms, respectively. Dashed lines represent the boundaries of the ``true" SBs. The energy threshold is $\xi_{\max}=0.1657$ and the external frequency is~$\Omega=0.89$\label{fig:saddle:connection}}
\end{figure}

Figure~\ref{fig:saddle:connection} shows an example of SB boundaries passing though the saddle connection for~$\xi_{\max}=0.1657$ and~$\Omega=0.89$. As one can see, when~$F<\widehat{F}$ there is a SBMT that wraps around the cylinder, hence, there exists a ``true" safe basin of maximum type. When~$F=\widehat{F}$, the maximum and the saddle curves coincide. When~$F > \widehat{F}$, SBMT becomes of kind I, hence, first, for a particular value of the phase~$\psi$ system exhibits the coexistence of SBMT and SBST, while the phase-invariant ``true" SB collapses.

\begin{figure}[H]
\centering
\subfloat[\label{fig:true:erosion:sim}]{\includegraphics[width=0.45\textwidth]{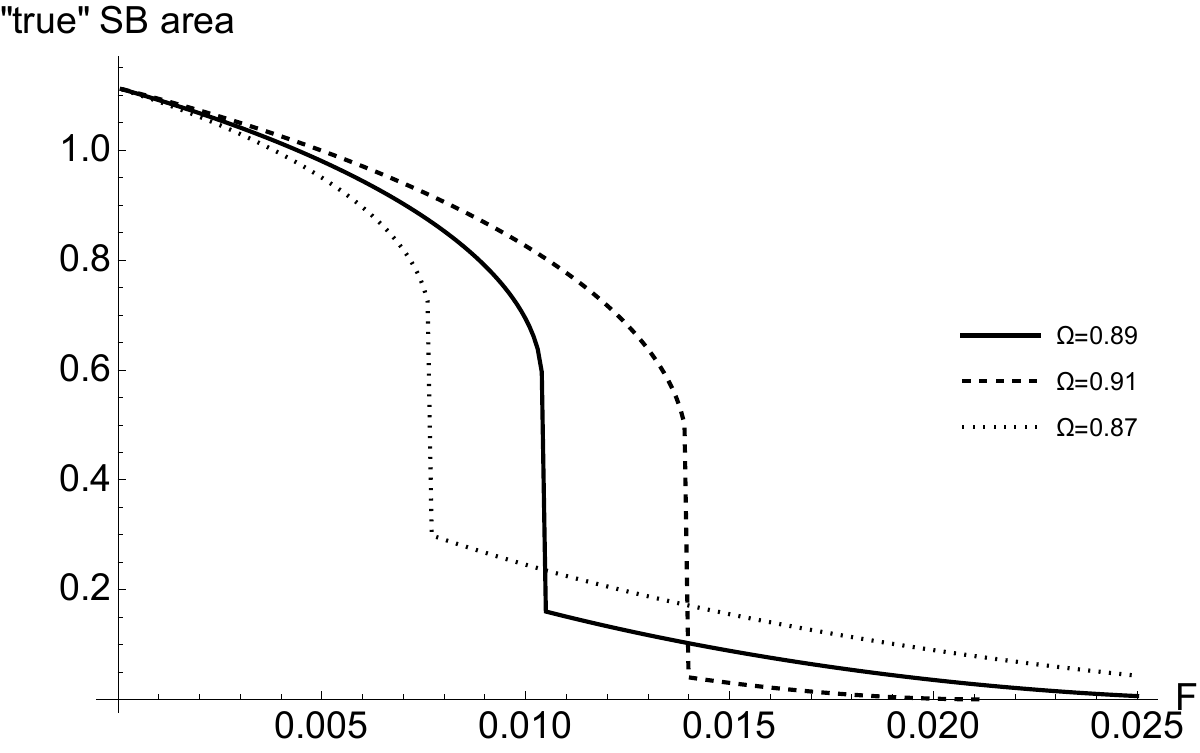}}
\hspace{10pt}
\subfloat[\label{fig:true:erosion:sim2}]{\includegraphics[width=0.45\textwidth]{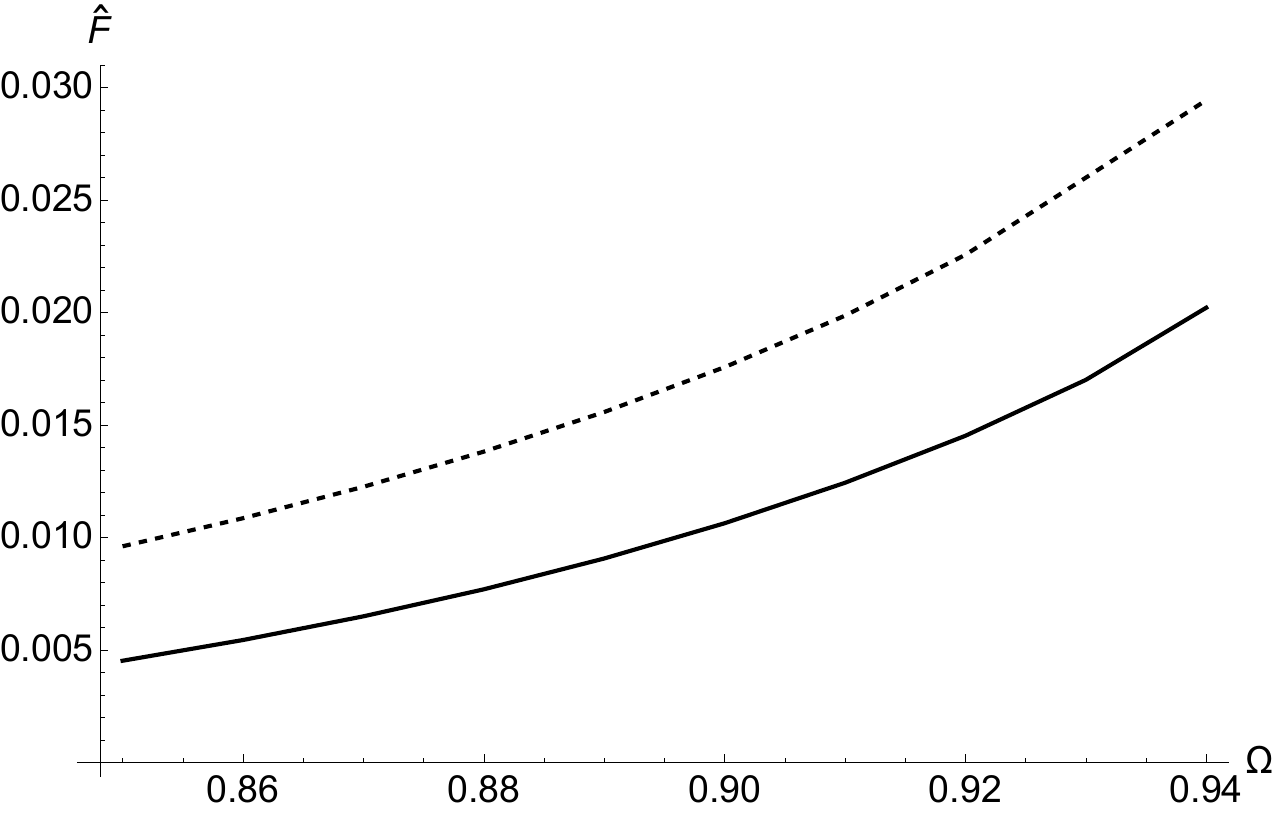}}
\caption{Panel (a): Approximation of the erosion profile for a truncation level $\xi=0.1567$. Panel (b): Critical value~$\widehat{F}$ as a function of the excitation frequency~$\Omega$. Dashed and solid curves correspond to~$\xi_{\max}=0.1657$ and~$\xi_{\max}=0.155$, respectively}
\end{figure}

Therefore, one can define an approximation of the ``true" SB erosion profile~$\mu(F)$, i.e., the GIM (area) of the ``true" SB as a function of the parameter~$F$, as
\[
\mu(F) = 2\pi J(\widehat{\xi})
\] where~$\widehat{\xi} = \widehat{\xi}(F)$ is the following piece-wise continuous function
\[
\widehat{\xi}(F) = \begin{cases} \xi_m(0), & F < \widetilde{F},\\  \xi_s(\pi), & F \ge \widetilde{F}.\end{cases}
\]
Figure~\ref{fig:true:erosion:sim} shows the erosion profiles of the ``true" SB, i.e., SB area depending on the increasing excitation frequency, for three values of the excitation frequency~$\Omega$. The truncation level is set to be~$\xi_{\max} = 0.158$. The critical forcing amplitude value~$\widehat{F}$ grows with increasing value of~$\Omega$ as demonstrated in Figure~\ref{fig:true:erosion:sim2}.

Unfortunately, the approximation of isolated $1:1$ resonance fails to capture highly irregular dynamics at the separatrix splitting, thus, rendering an erroneous approximation near the top of the RM cylinder. In order to visualize the discrepancy between the actual escaping initial conditions and the RM, one can, first, discretize the~$(\vartheta, \xi)$-cylinder into a two-dimensional grid. Then, for each point on the grid perform the transformation~\eqref{eq:q:trans}--\eqref{eq:p:trans}, and use the obtained ICs to simulate system~\eqref{eq:main} until either the escape occurs or the simulation time exceeds a preset maximum time~$T_{\max}$. Assigning a color to each point according to whether escape occurred or not, allows one to obtain a color plot of the actual escaping ICs on the~$(\vartheta, \xi)$-cylinder.

\begin{figure}[H]
\centering
\subfloat[$F=0.008$]{\includegraphics*[width=0.3\textwidth]{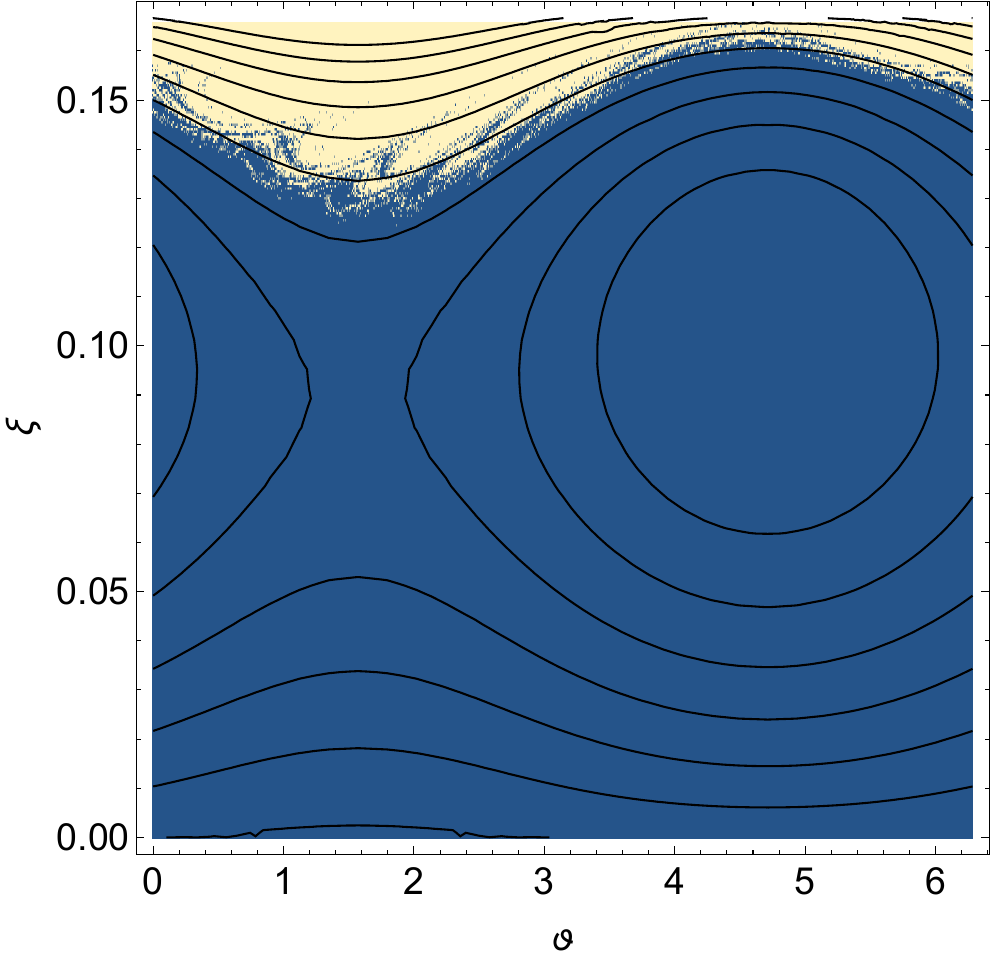}}
\hspace{10pt}
\subfloat[$F=0.012$]{\includegraphics*[width=0.3\textwidth]{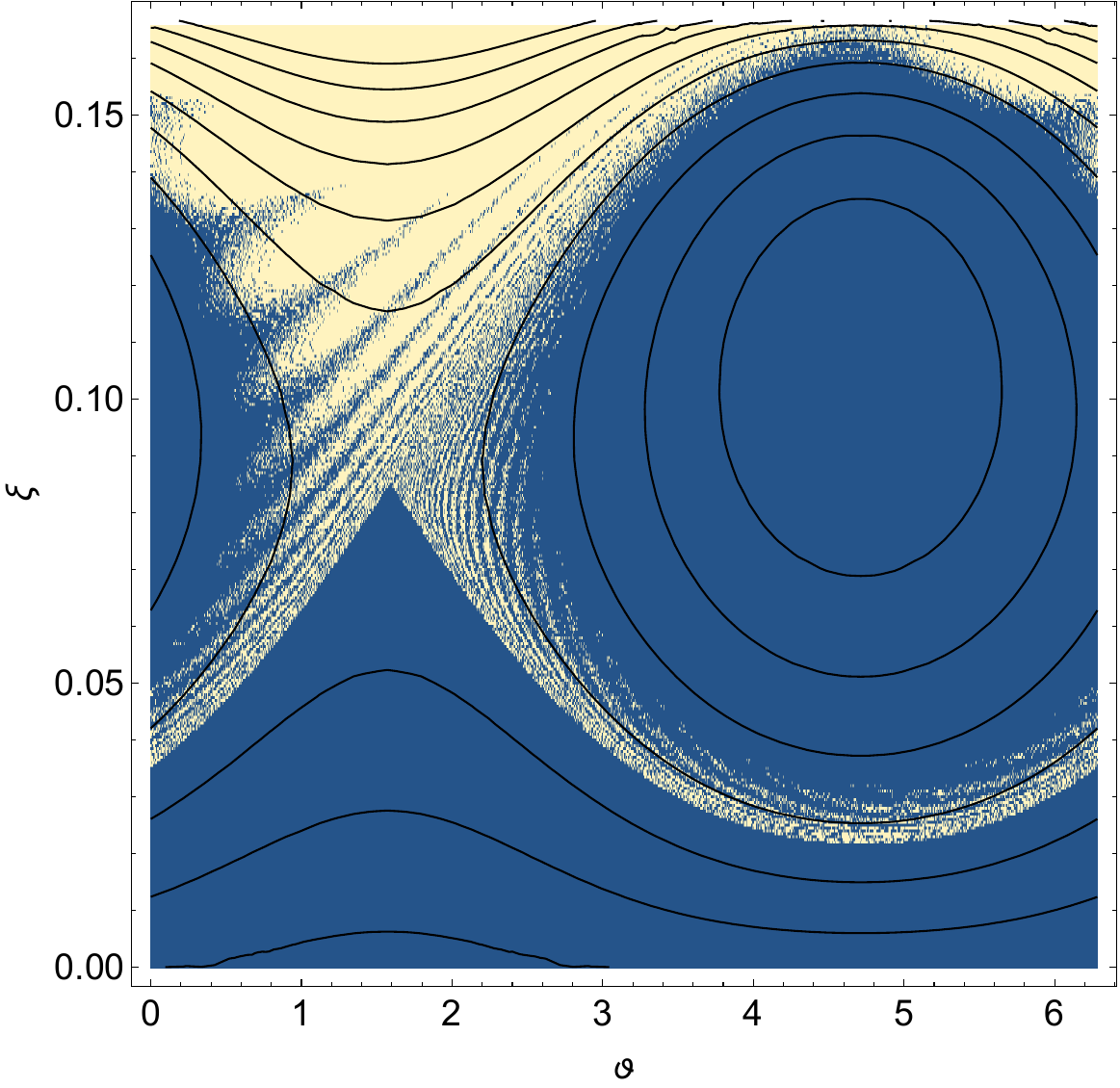}}
\hspace{10pt}
\subfloat[$F=0.016$]{\includegraphics*[width=0.3\textwidth]{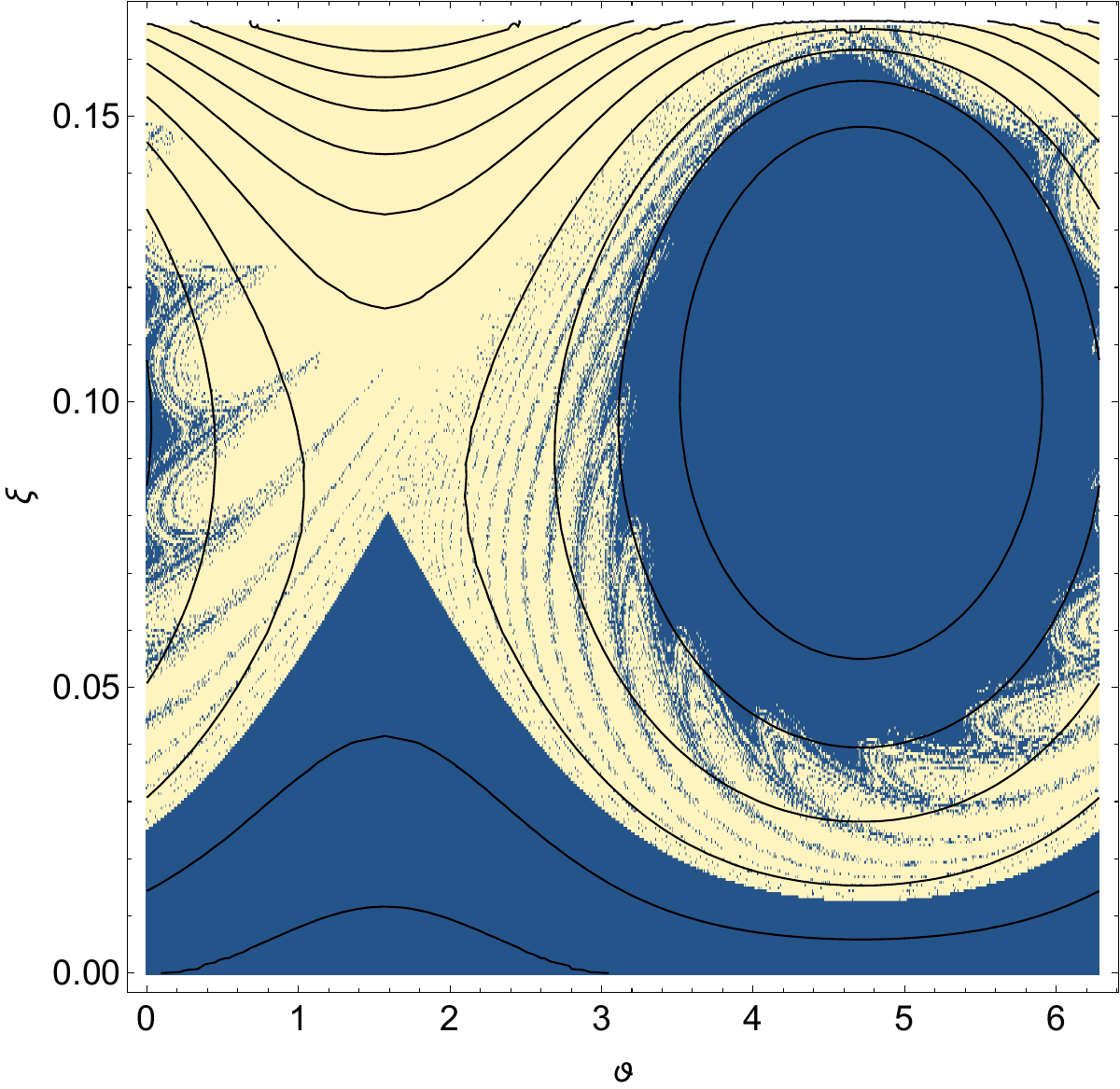}}
\caption{Actual escaping ICs (yellow) plotted on the $(\vartheta, \xi)$-cylinder\label{fig:discrepancy:cylinder}}
\end{figure}

Figure~\ref{fig:discrepancy:cylinder} presents the actual escaping initial conditions (yellow) overlayed with the resonance manifold for three different values of forcing amplitude~$F$. Visualizing actual escape as a color plot on the RM yields several critical insights. First of all, even for a small value of~$F$ the maximum mechanism on the RM fails to capture the actual escape dynamics, and therefore, some additional steps are required in order to account for this effect. Secondly, as~$F$ increases, we see the progression of the erosion of SBMT while SBST stays intact with an accurate robust approximation by the slow-flow motion. Finally, Figure~\ref{fig:discrepancy:cylinder} confirms that despite the mentioned discrepancy, the $(\vartheta, \xi)$-variables are well suited for the study of the SB. In particular, it is evident that the erosion of the SB can be explained through the dynamics on the RM as it undergoes the separatrix splitting. Note that, the RM on Figure~\ref{fig:discrepancy:cylinder} is shifted horizontally by $\pi/2$.

To overcome the limitation of the maximum mechanism, one can utilize the truncation of the potential, i.e., using a lower energy level as the escape threshold. In order to find the appropriate truncation level one can search for a limiting level curve of~$C(\vartheta, \xi)$ around the center critical point such that all the points on the curve are non-escaping (after applying the transformation~\eqref{eq:q:trans}--\eqref{eq:p:trans}). Taking the maximum value~$\xi^*$ of $\xi$ on this curve provides the {\em effective}, i.e., corrected, energy threshold value to be used for the potential truncation. Note that, the new boundary of SBMT on the RM becomes:
\[
C(\vartheta, \xi) = C(\vartheta^*, \xi^*),
\] where~$\vartheta = \pi$.
In order to numerically obtain the effective threshold value~$\xi^*$ one can adapt a bisection method based on the following observations.

First, notice that there is a critical point~$\left(\vartheta^*, \xi_c\right)$ of the center type with the phase curves of the slow-flow system ``wrapped" around it. If the point~$\left(\vartheta^*, \xi_c\right)$ is the IC of a non-escaping trajectory (in the full system), so are all the points in a sufficiently small neighborhood of it. Secondly, one can choose a point~$\left(\vartheta^*, \xi^\dagger\right)$ close to~$\xi=\xi_{\max}$ that corresponds to the escaping trajectory. Since the value of the slow-phase~$\vartheta$ is the same, consider an interval~$\mathcal{I} = [\xi_c,\; \xi^\dagger]$. For each $\widetilde{\xi} \in \mathcal{I}$ one can derive the phase curve
\[
\Gamma_{\widetilde{\xi}} = \left\{(\vartheta, \xi)\; : \;C(\vartheta, \xi) = C(\vartheta^*, \widetilde{\xi})\right\},
\] and perform the transformation~\eqref{eq:q:trans}--\eqref{eq:p:trans} to locate the corresponding curve~$\Gamma_{pq}$ on the~$(q,p)$-plane. Then, one can define a boolean function~$\Upsilon(q,p) : \mathbb{R}^2 \to \left\{0,1\right\}$ according to whether the trajectory starting at~$(q, p)$ is escaping or not, and apply it to every point of~$\Gamma_{pq}$. Thus, one can construct a mapping~$\mathfrak{M}: \mathcal{I} \to \{0,1\}$ as follows:
\[
\mathfrak{M}\left(\widetilde{\xi}\right) = \bigwedge\limits_{(\vartheta, \xi)\in\Gamma_{\widetilde{\xi}}} \Upsilon\left(q(\vartheta, \xi),\; p(\vartheta, \xi)\right),
\] where symbol~$\wedge$ denotes the conjunction, i.e., the logical operator~``\texttt{and}". Since~$\mathfrak{M}$ maps the interval~$\mathcal{I}$ to a discrete topological space~$\{0,1\}$, and mapping $\xi \mapsto C(\vartheta^*, \xi)$ is continuous, $\mathfrak{M}$ is locally constant, and therefore, also continuous. Moreover, it is easy to see that function~$\mathfrak{M}$ takes the opposite values on the endpoints of the interval~$\mathcal{I}$, i.e., $\mathfrak{M}(\xi_c) = 0$ and~$\mathfrak{M}(\xi^\dagger) = 1$. Thus, one can apply the bisection algorithm with some preset tolerance~$\epsilon$ to obtain the limiting value~$\xi^*$ such that $\mathfrak{M}(\xi) = 0$ for every $\xi \le \xi^*$ and $\mathfrak{M}(\xi) = 1$ for all~$\xi \ge \xi^* + \epsilon$. In other words, all the points on the curve~$C\left(\vartheta, \xi\right) = C(\vartheta^*, \xi^*)$ correspond to the non-escaping trajectories in the full system, while for any~$\xi\ge \xi^* + \epsilon$ the corresponding phase curve contains at least one initial condition that leads to escape. Therefore, we obtain the corrected value of the escape threshold~$\xi = \xi^*$.

Apart from setting the tolerance~$\epsilon$, one needs to select the maximal simulation time~$T_{\max}$ and the grid size~$\delta$ for the discretization of the phase curves at each iteration of the method. Choosing an appropriate value of~$\delta$ is crucial for rendering reliable results. Clearly, the accuracy of the algorithm improves with decreasing~$\delta$. However, calculating~$\mathfrak{M}$ is computationally demanding, and opting for a smaller~$\delta$ can substantially impede the computational speed.

\begin{figure}[H]
\centering
\subfloat[\label{fig:trunc:levels}]{\includegraphics*[width=0.48\columnwidth]{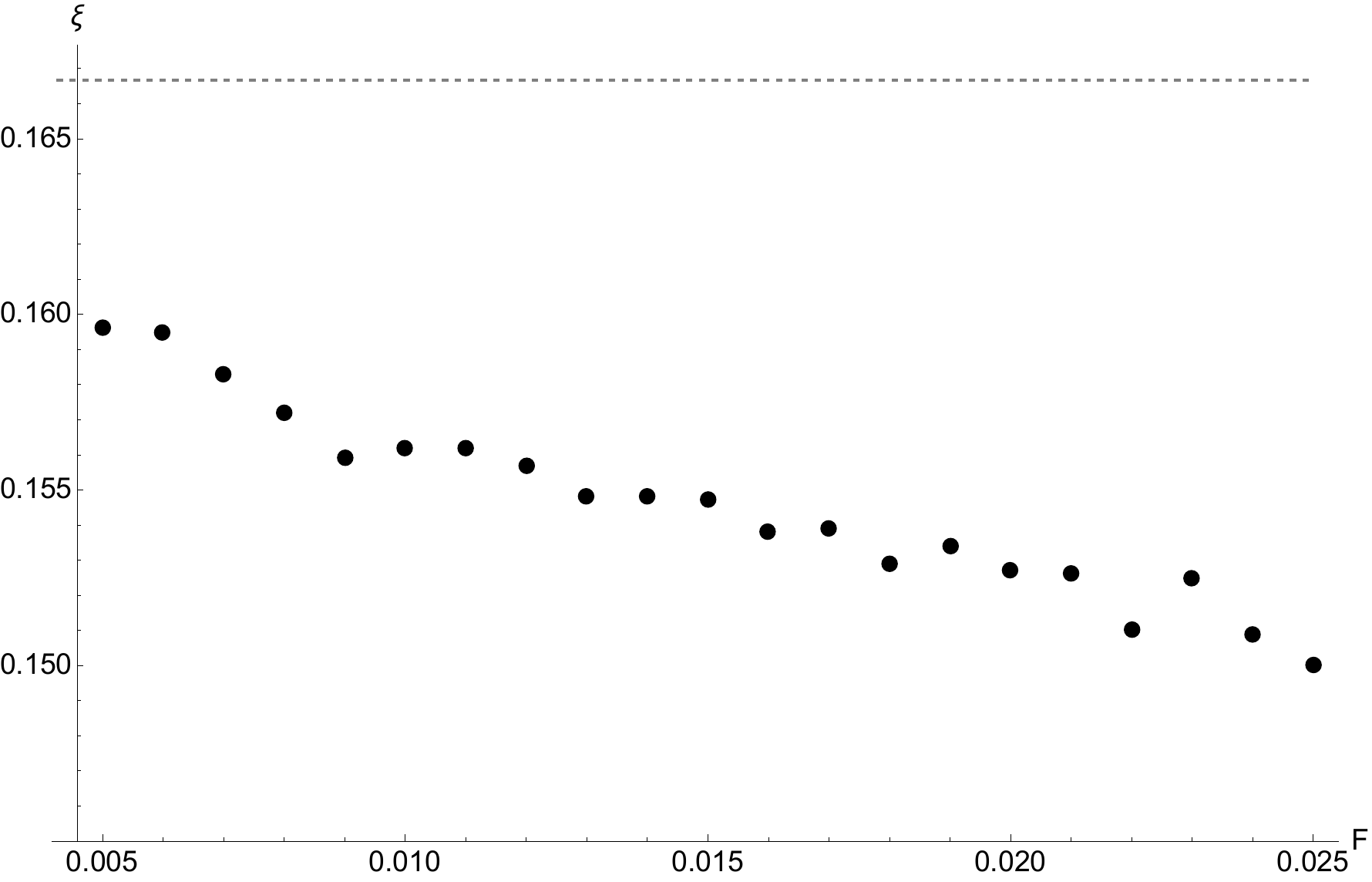}}
\hspace{10pt}
\subfloat[\label{fig:area:comparison:phase:0}]{\includegraphics[width=0.48\columnwidth]{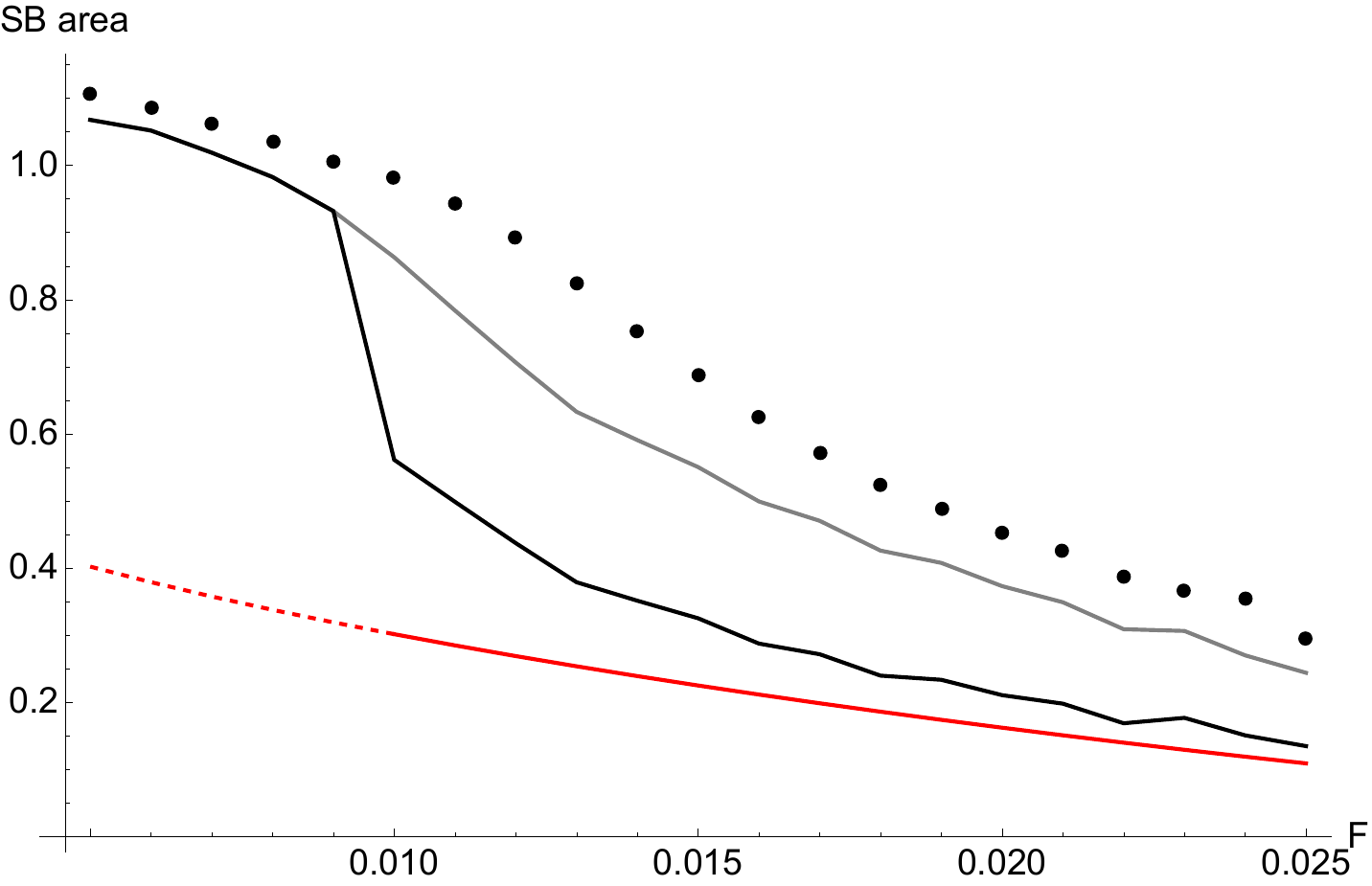}}
\caption{Panel (a): Black dots represent numerically obtained thresholds that account for the stochastic layer. Gray dashed line is the maximum depth of the potential well, i.e.,~$\xi=1/6$. Panel (b): Erosion profiles, i.e., SB areas depending on the forcing amplitude $F$. Black curve represents the area of SBMT. Red curve is the area of SBST (dashed part denotes area of SBST which is completely included inside SBMT). Black dots correspond to the area of SB computed numerically. Grey curve represents the sum of the areas of the SB of both types.}
\end{figure}

\begin{figure}[hbtp]
\centering
\subfloat[$F=0.007$]{\includegraphics*[width=0.3\textwidth]{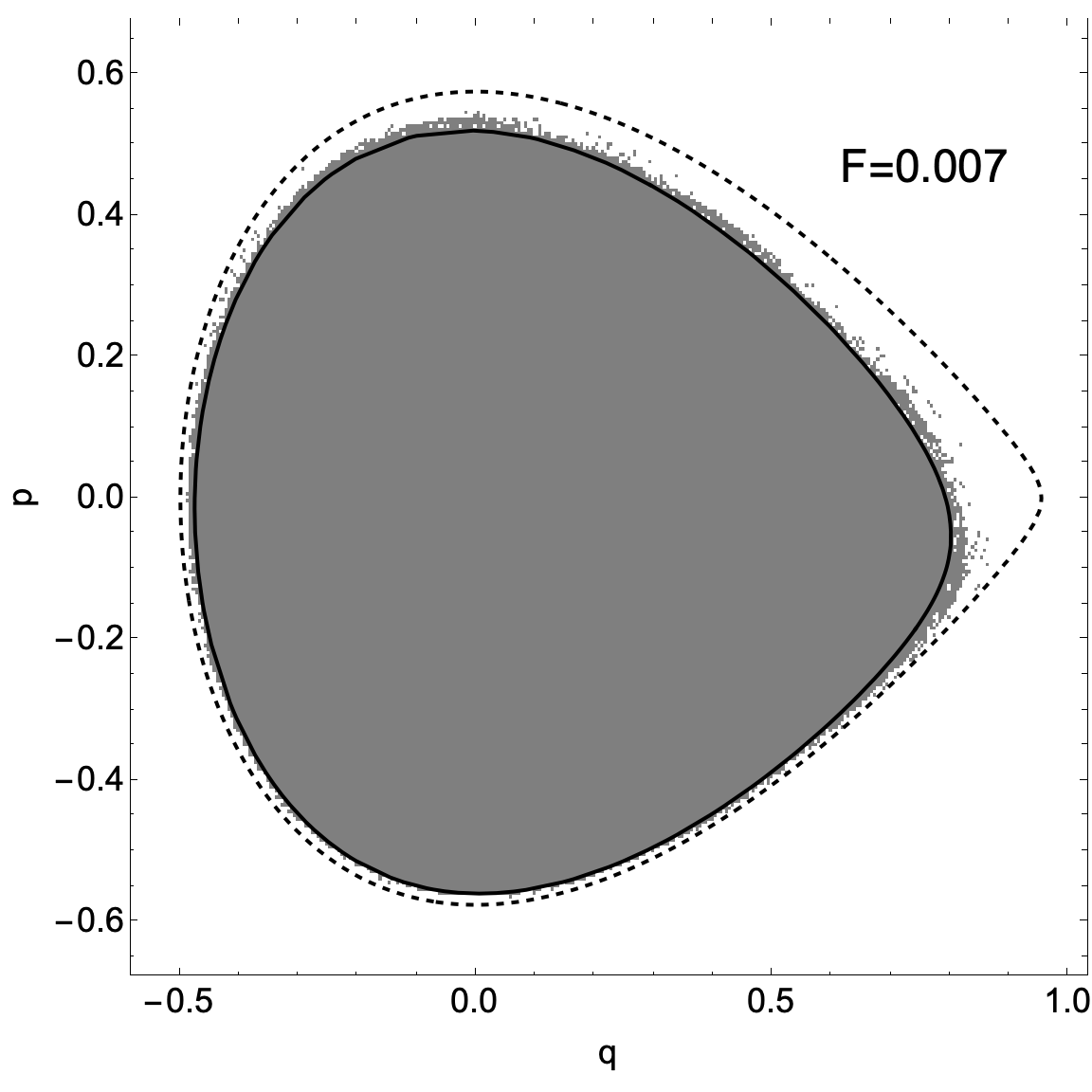}}
\hspace{10pt}
\subfloat[$F=0.009$]{\includegraphics*[width=0.3\textwidth]{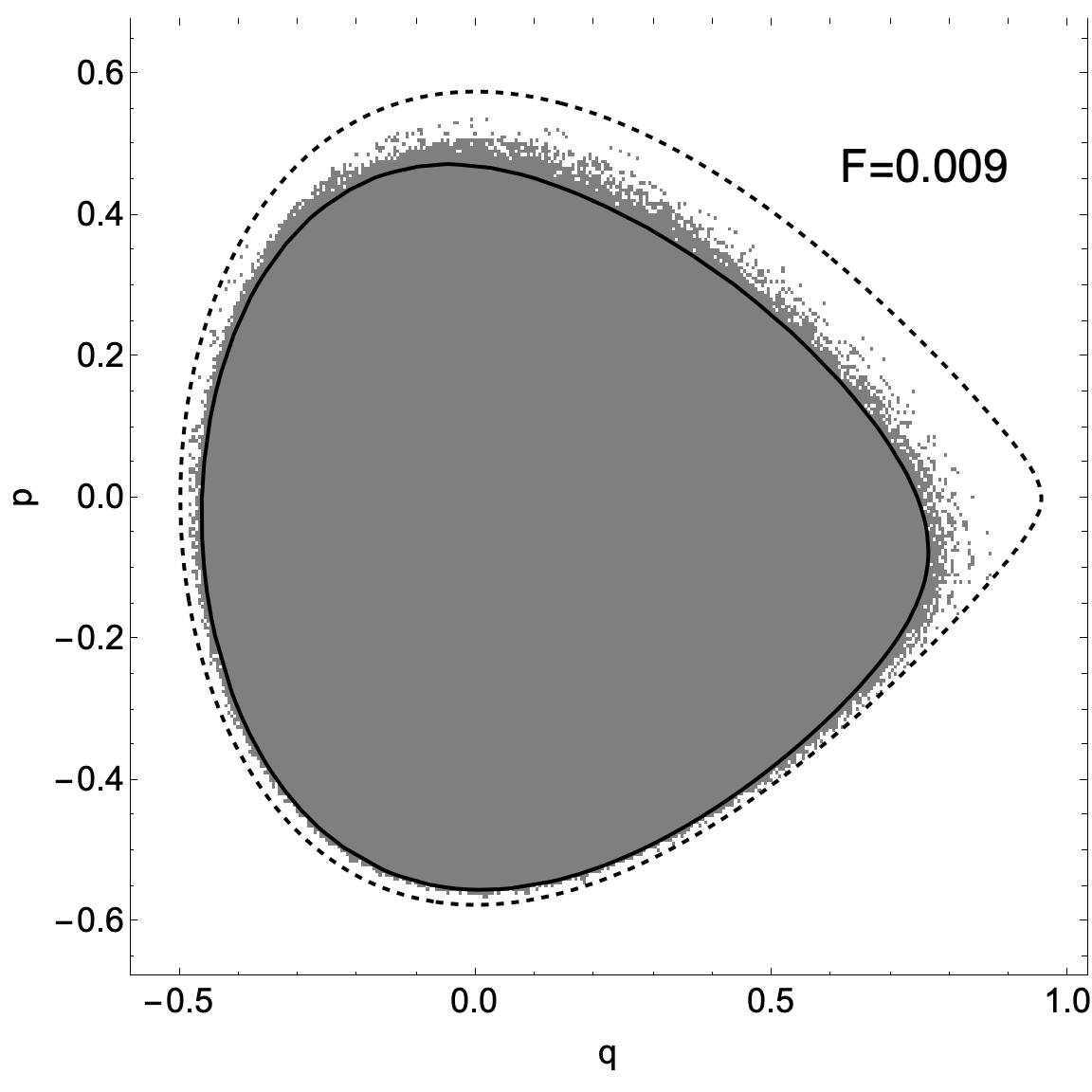}}
\hspace{10pt}
\subfloat[$F=0.01$]{\includegraphics*[width=0.3\textwidth]{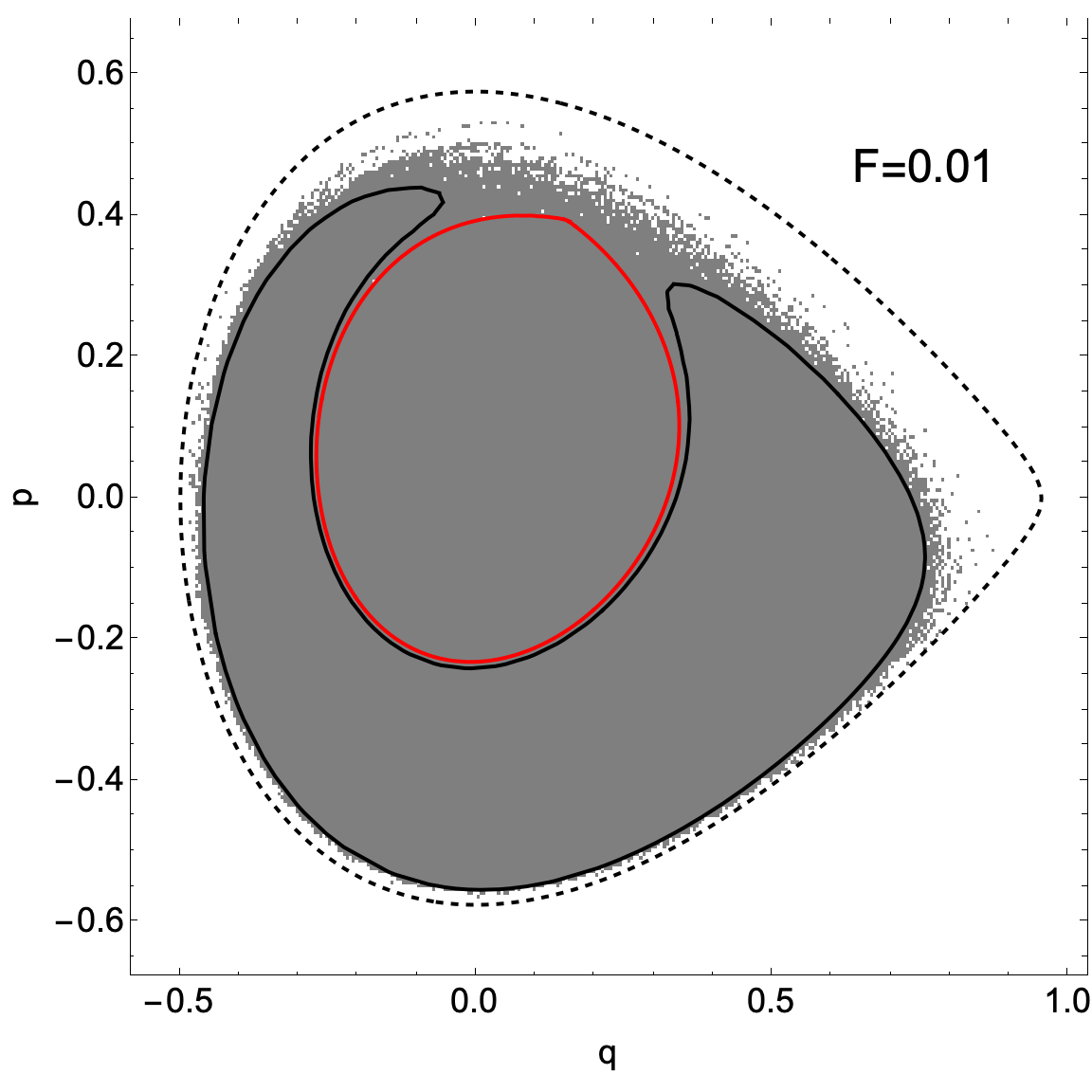}}\\
\subfloat[$F=0.011$]{\includegraphics*[width=0.3\textwidth]{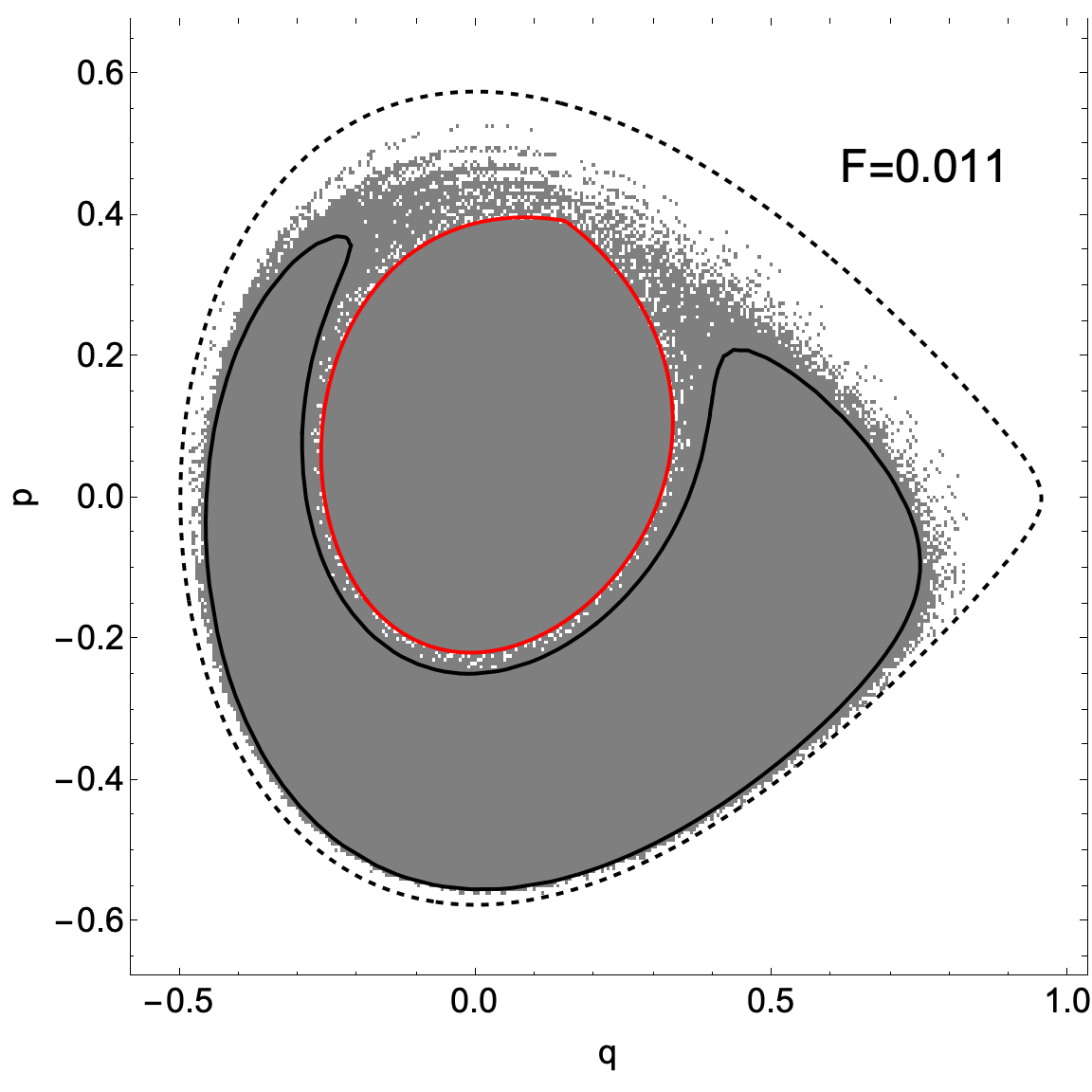}}
\hspace{10pt}
\subfloat[$F=0.012$]{\includegraphics*[width=0.3\textwidth]{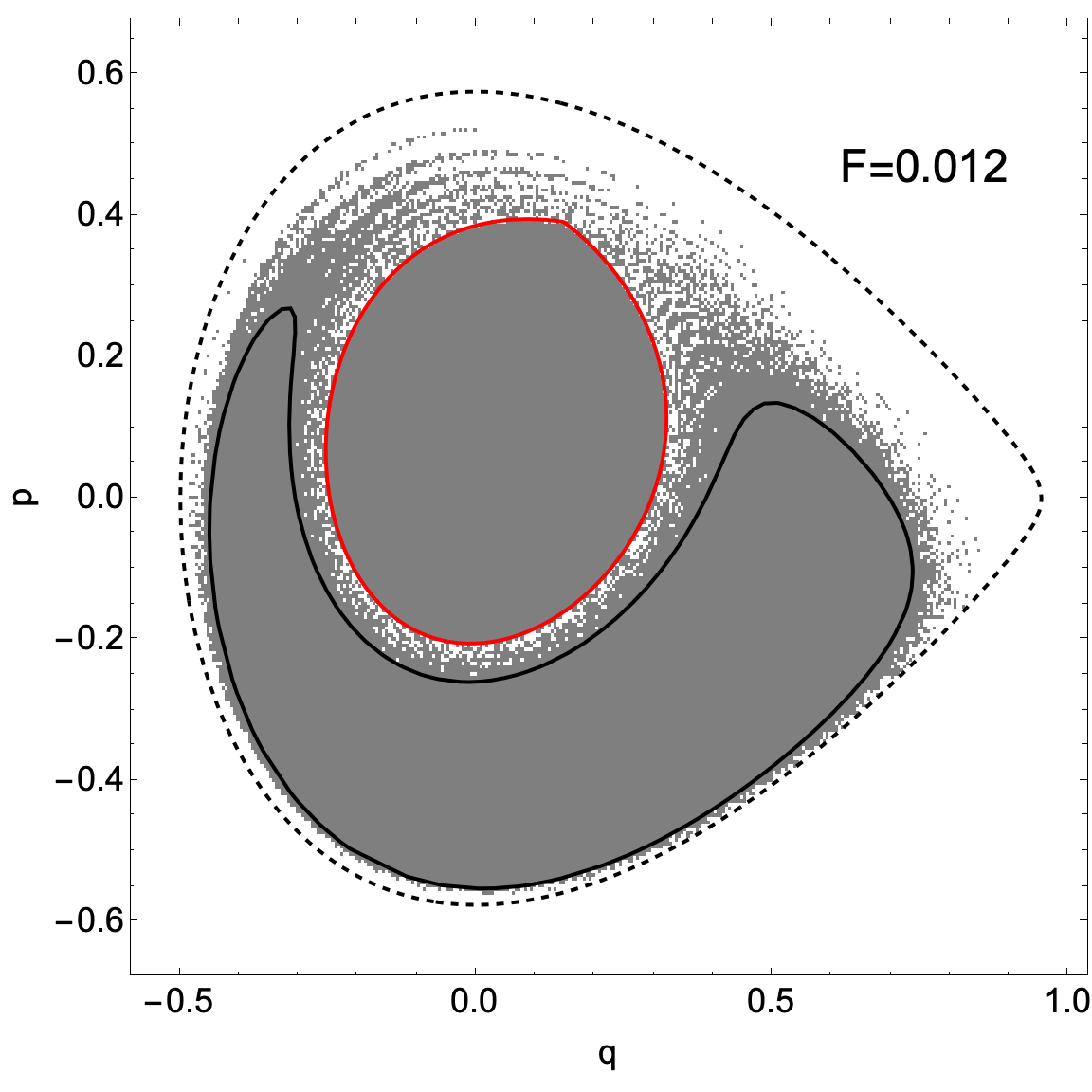}}
\hspace{10pt}
\subfloat[$F=0.013$]{\includegraphics*[width=0.3\textwidth]{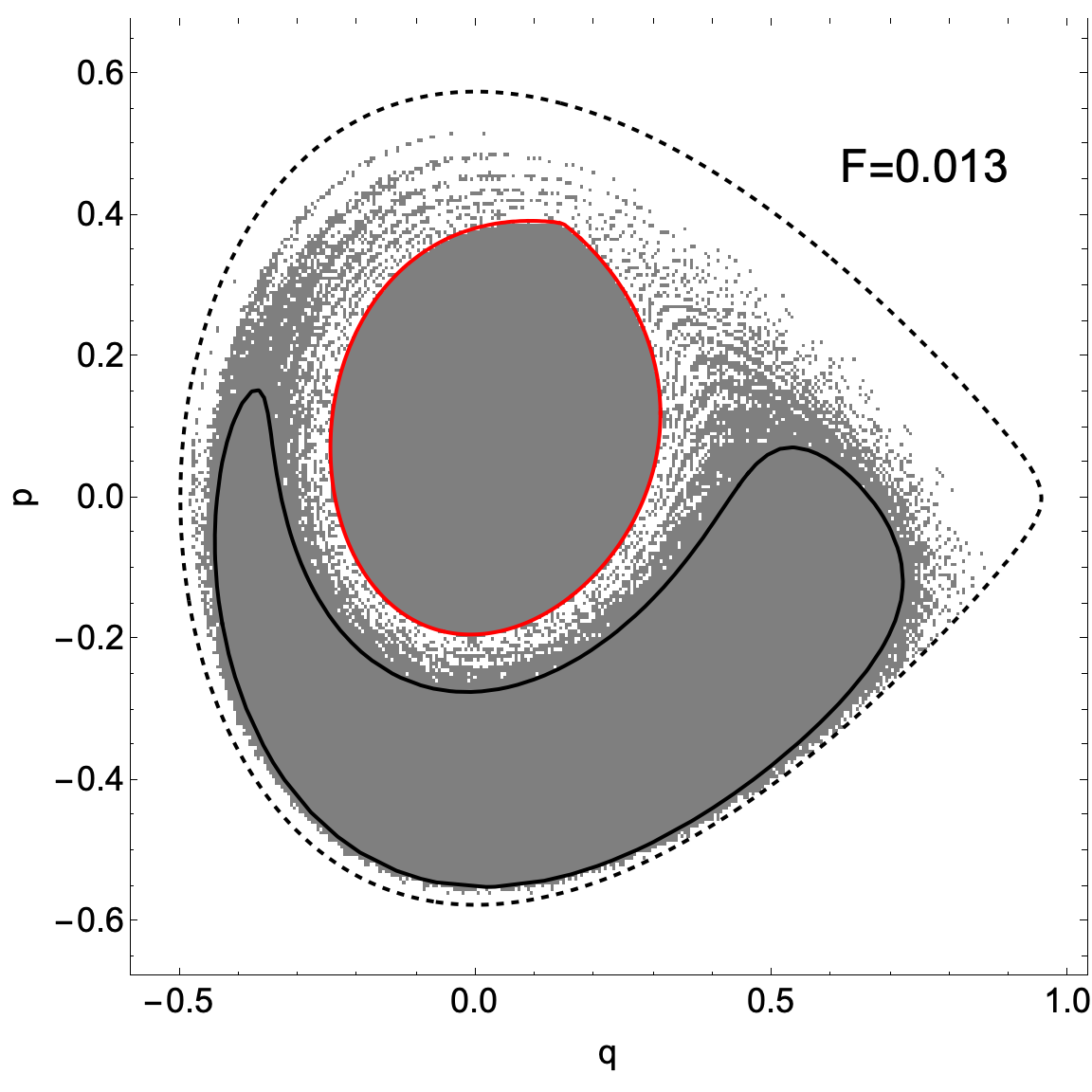}}\\
\subfloat[$F=0.016$]{\includegraphics*[width=0.3\textwidth]{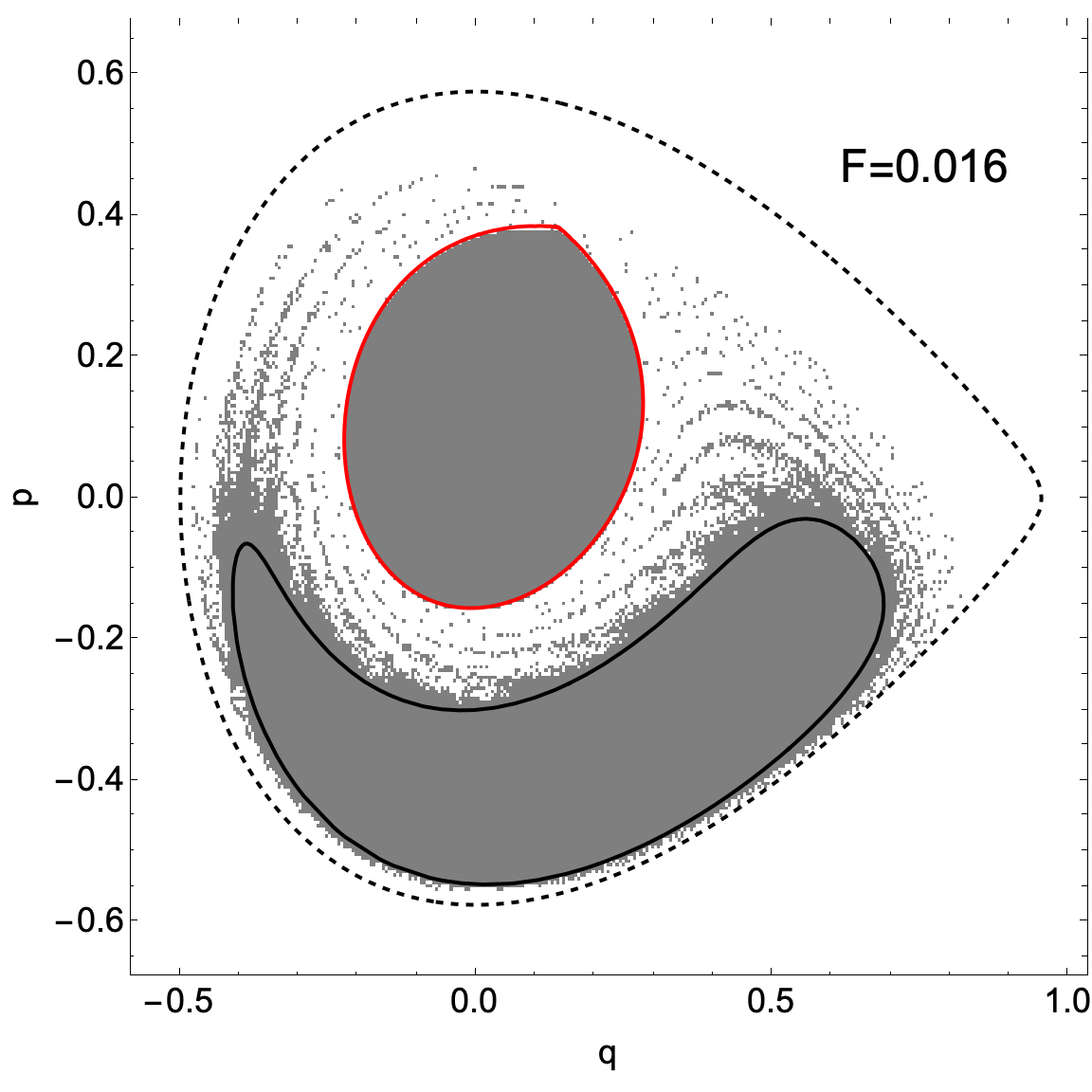}}
\hspace{10pt}
\subfloat[$F=0.02$]{\includegraphics*[width=0.3\textwidth]{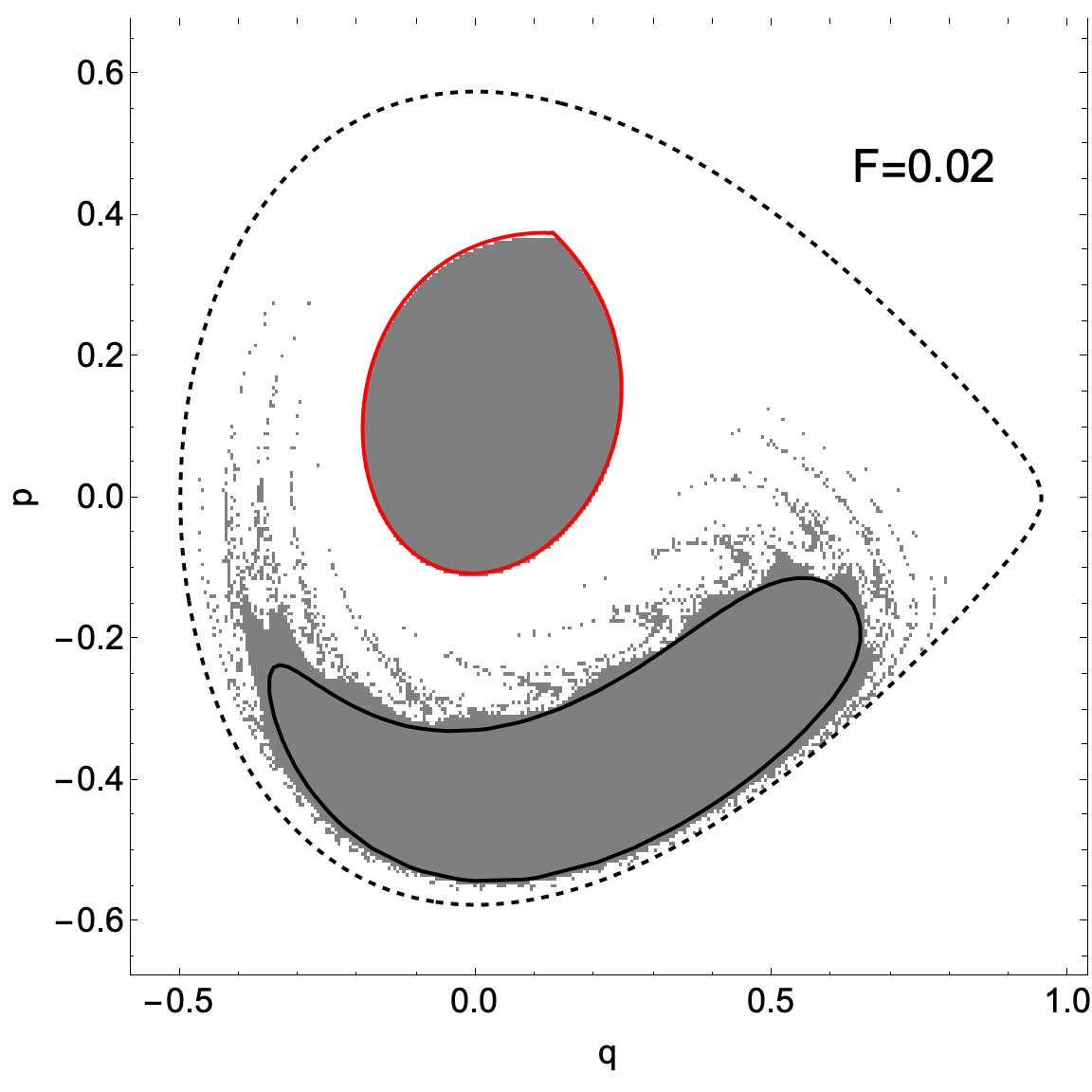}}
\hspace{10pt}
\subfloat[$F=0.025$]{\includegraphics*[width=0.3\textwidth]{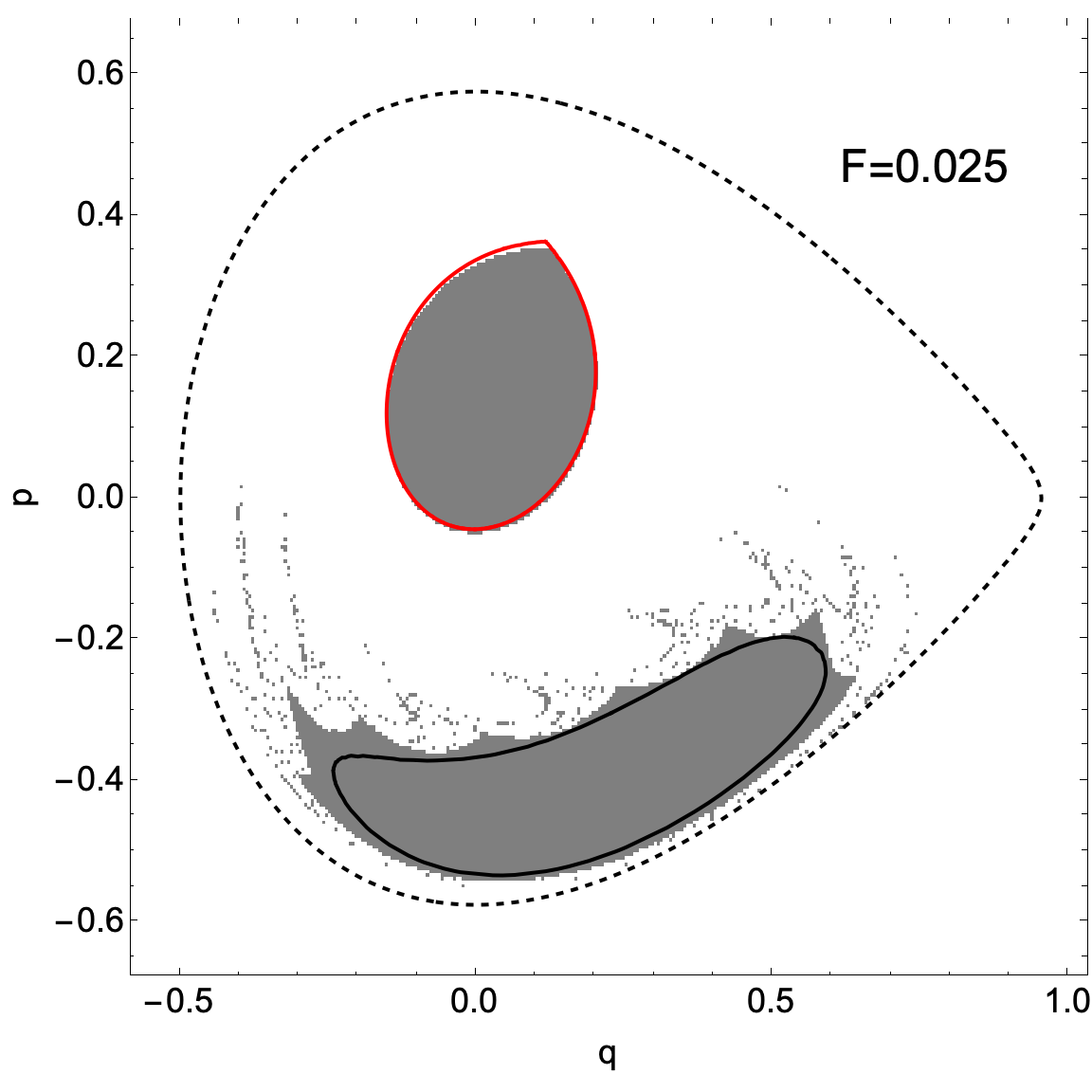}}
\caption{Evolution of the safe basin for the phase $\psi=0$. Gray color represents the initial conditions for which the numerical simulations show no escape. Each numerical simulation runs for $100$ periods of the external excitation with the frequency~$\Omega=0.89$\label{fig:sb:evolution:phase:0}}
\end{figure}

Figure~\ref{fig:trunc:levels} shows numerically obtained threshold values depending on the external forcing amplitude~$F$. For the computations the following parameters were used: tolerance~$\epsilon = 0.001$, $T_{\max} = 100$~periods, and~$\delta = 0.01$. As one can see, the discrepancy between the effective threshold values and the maximum depth of the potential well increases with~$F$. Using obtained values allows one to account for the limitation of the maximum escape mechanism and render a better approximation of the SB boundaries, see Figure~\ref{fig:sb:evolution:phase:0}.

Figure~\ref{fig:sb:evolution:phase:0} presents an evolution of both numeric and approximate safe basins for the phase $\psi=0$. The boundaries of SBMT and SBST are presented with a black and red curves, respectively. As demonstrated, the approximate SBs yields a superb correspondence to the numerics by capturing the exact shape of the SB while leaving out the erosion. Note that, around $F=0.01$, the SB splits into two connected components and the AIR method allows us to precisely locate both SBs even when the split is not apparent from numerical simulations.

Furthermore, Figure~\ref{fig:area:comparison:phase:0} shows a comparison of the total area of the approximate and the numerically obtained SBs. As one can observe, the shape of the approximate erosion profile matches the numerics. However, all the numerical calculations of the area yield greater values since they include the fractal part of the SB which in fact is not considered safe.

\begin{figure}[H]
\centering
\subfloat[\label{fig:gim:true:comparison}]{\includegraphics*[width=0.5\columnwidth]{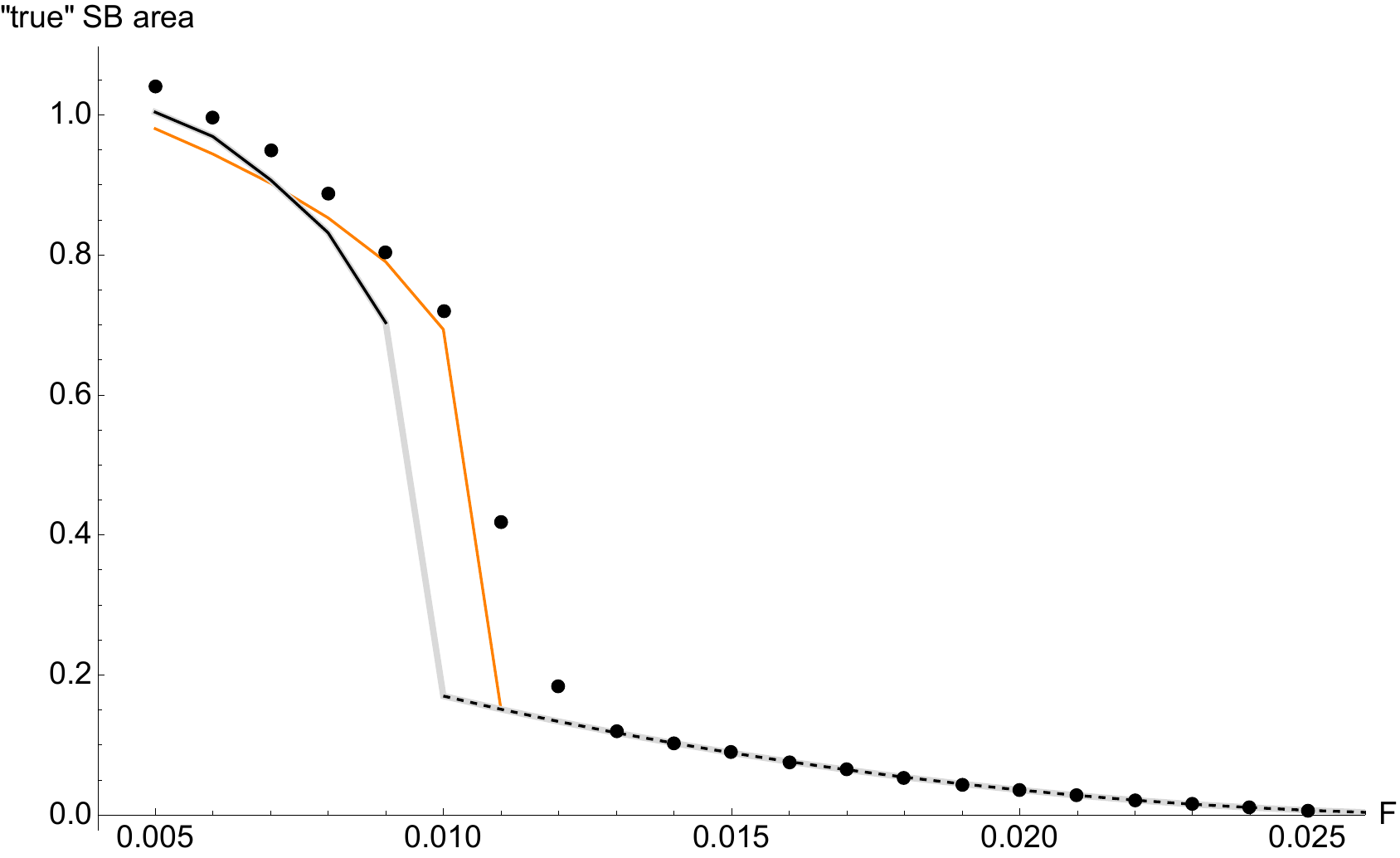}}
\hspace{10pt}
\subfloat[\label{fig:gim:true:comparison:zoom}]{\includegraphics*[width=0.4\columnwidth]{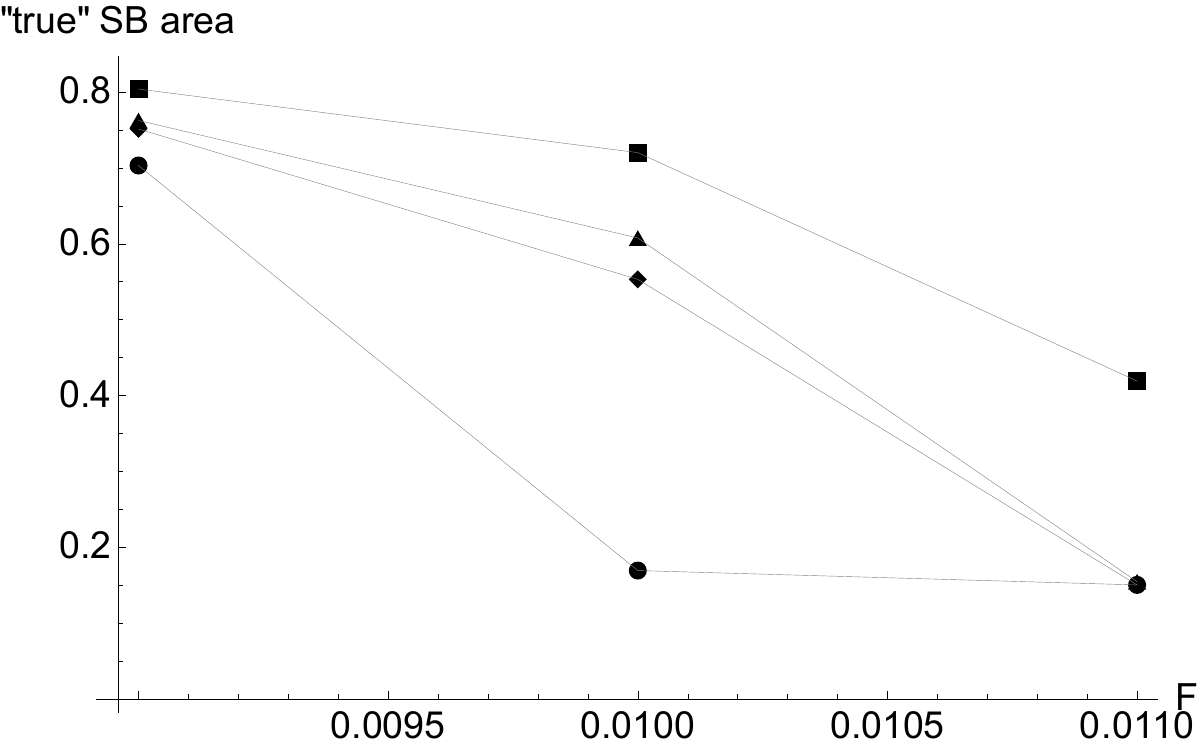}}
\caption{Comparison of the numerical and approximate erosion profiles of the ``true" safe basin. Panel (a): Black dots correspond to the numerical values of the ``true" SB areas (taken as an intersection of 21 steps along the period). Orange curve is an approximate erosion profile using energy threshold value~$\xi_{\max} = 0.158$. Panel (b): squares, triangles, diamonds correspond to numerical values of the ``true" SB areas taken for 20, 100 and 200 steps along the phase~$\psi$, respectively; circles correspond to the approximate values of the ``true" SB areas.\label{fig:true:sb:area:comparison}}
\end{figure}

\begin{figure}[H]
\centering
\subfloat[$F=0.007$]{\includegraphics*[width=0.3\textwidth]{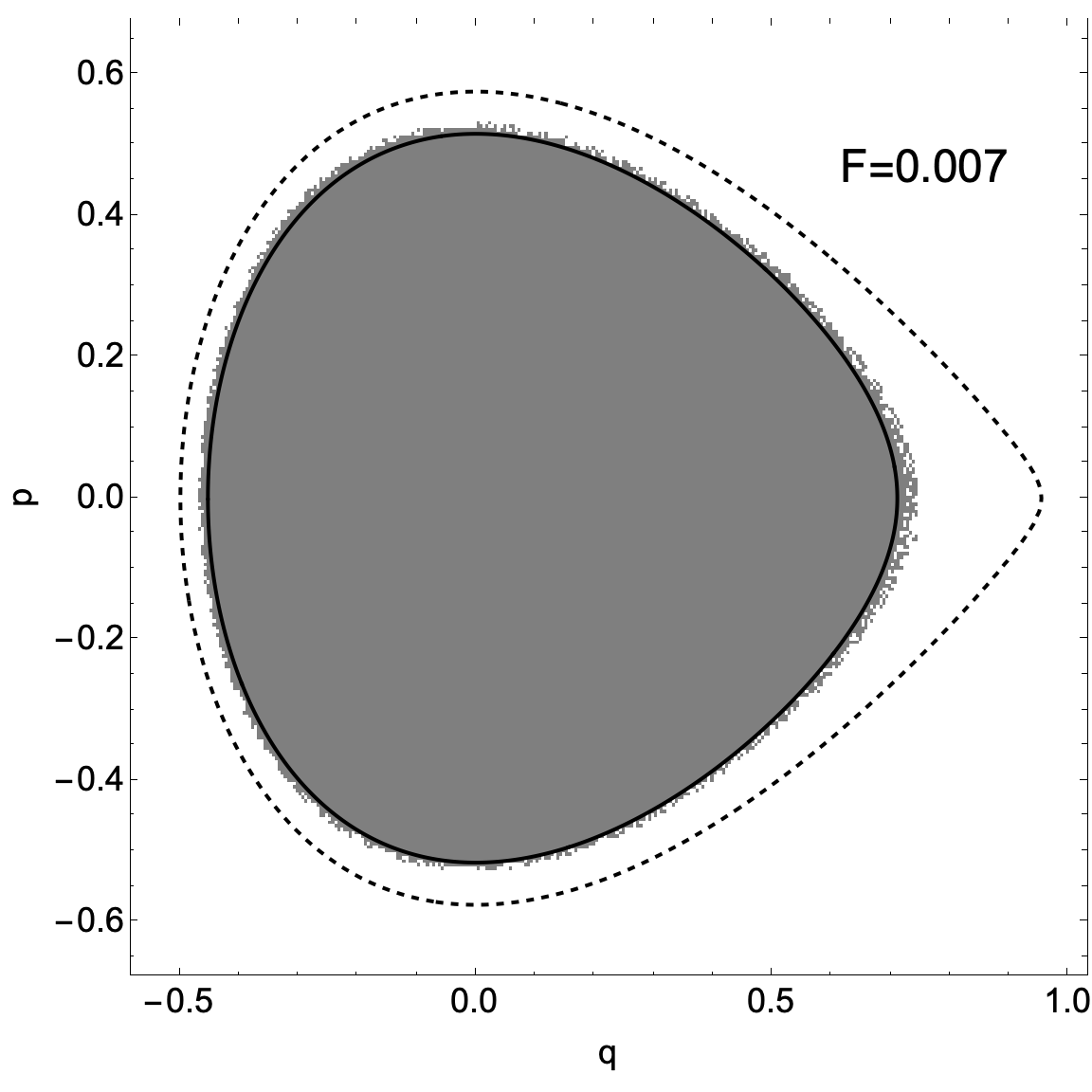}}
\hspace{10pt}
\subfloat[$F=0.009$]{\includegraphics*[width=0.3\textwidth]{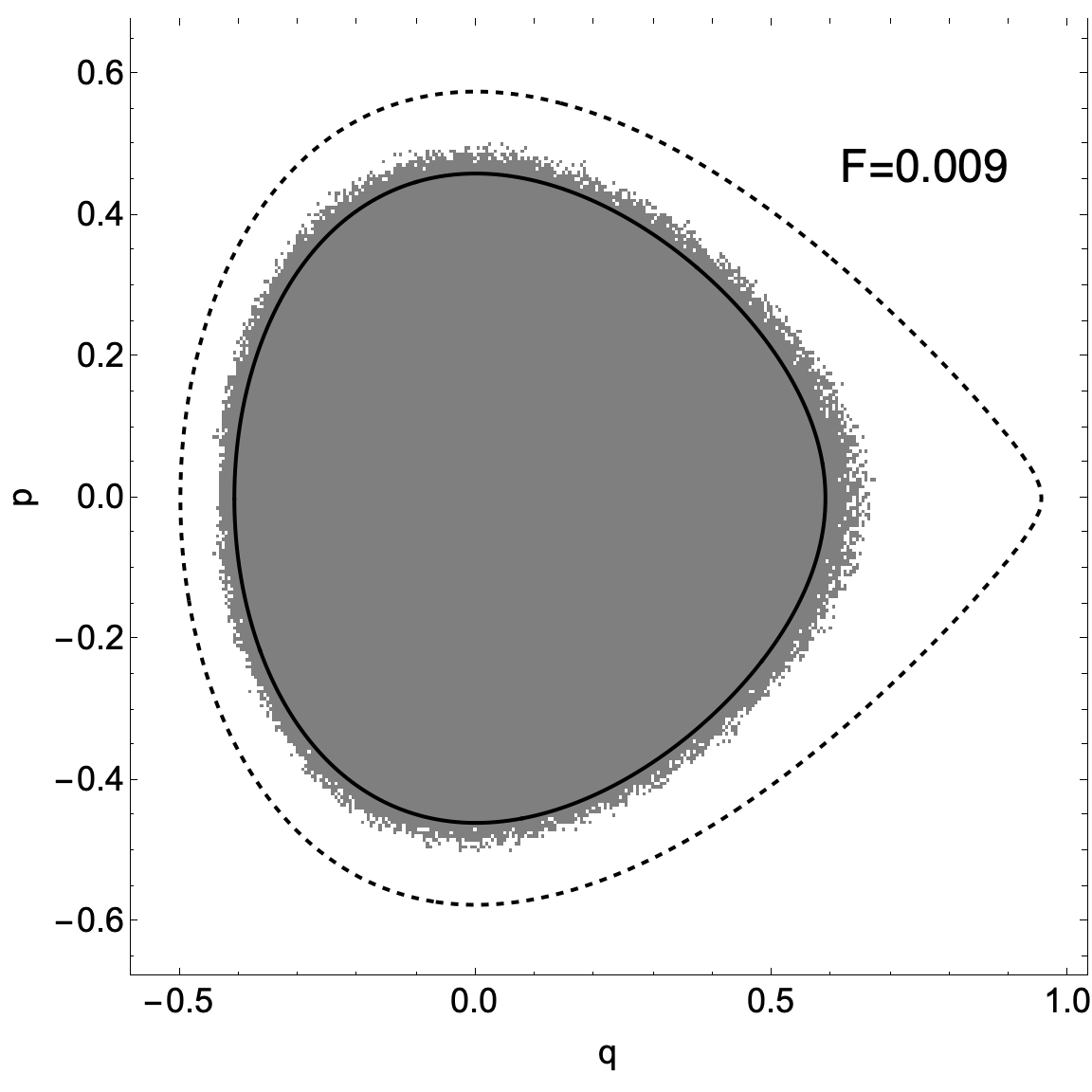}}
\hspace{10pt}
\subfloat[$F=0.01$]{\includegraphics*[width=0.3\textwidth]{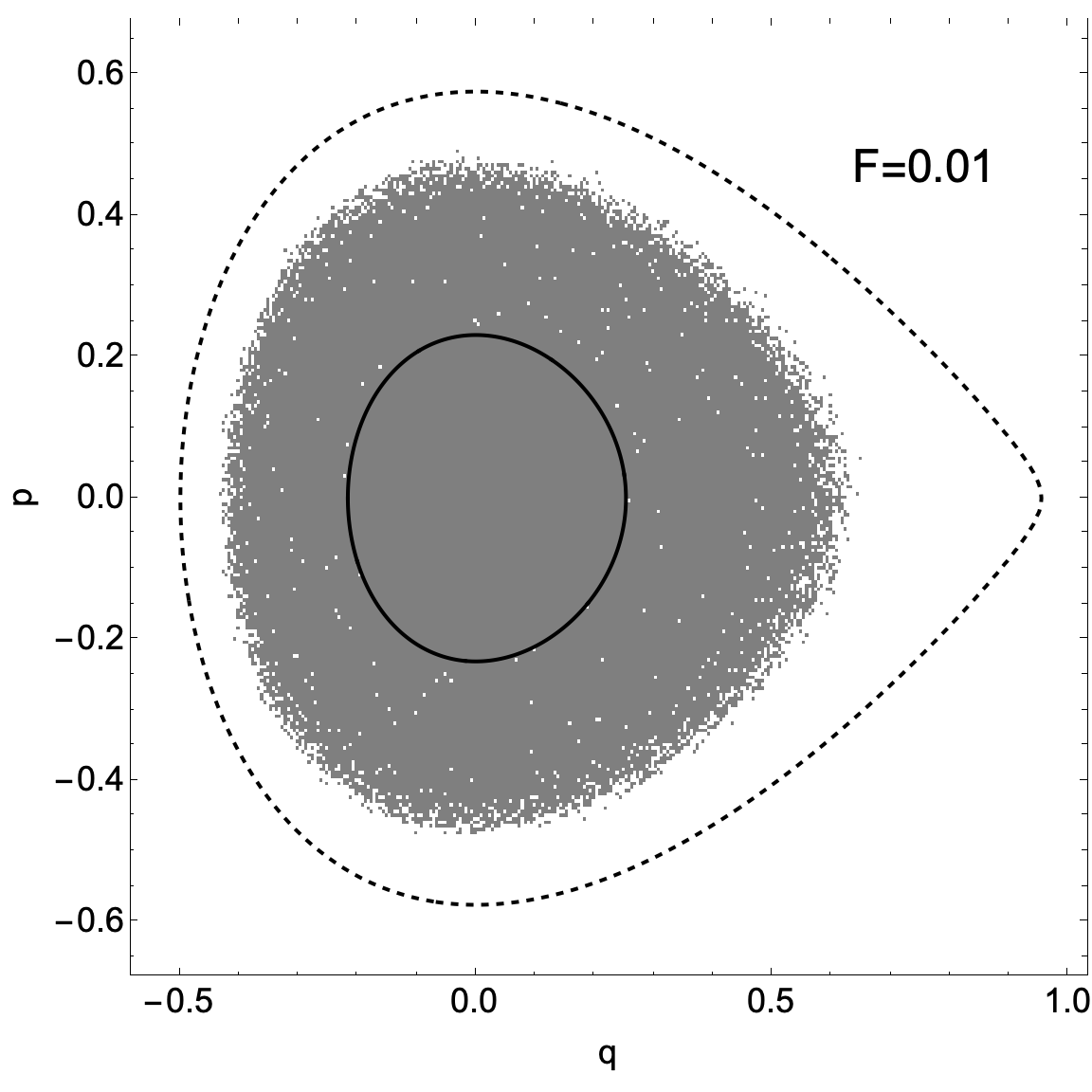}}\\
\subfloat[]{\includegraphics*[width=0.3\textwidth]{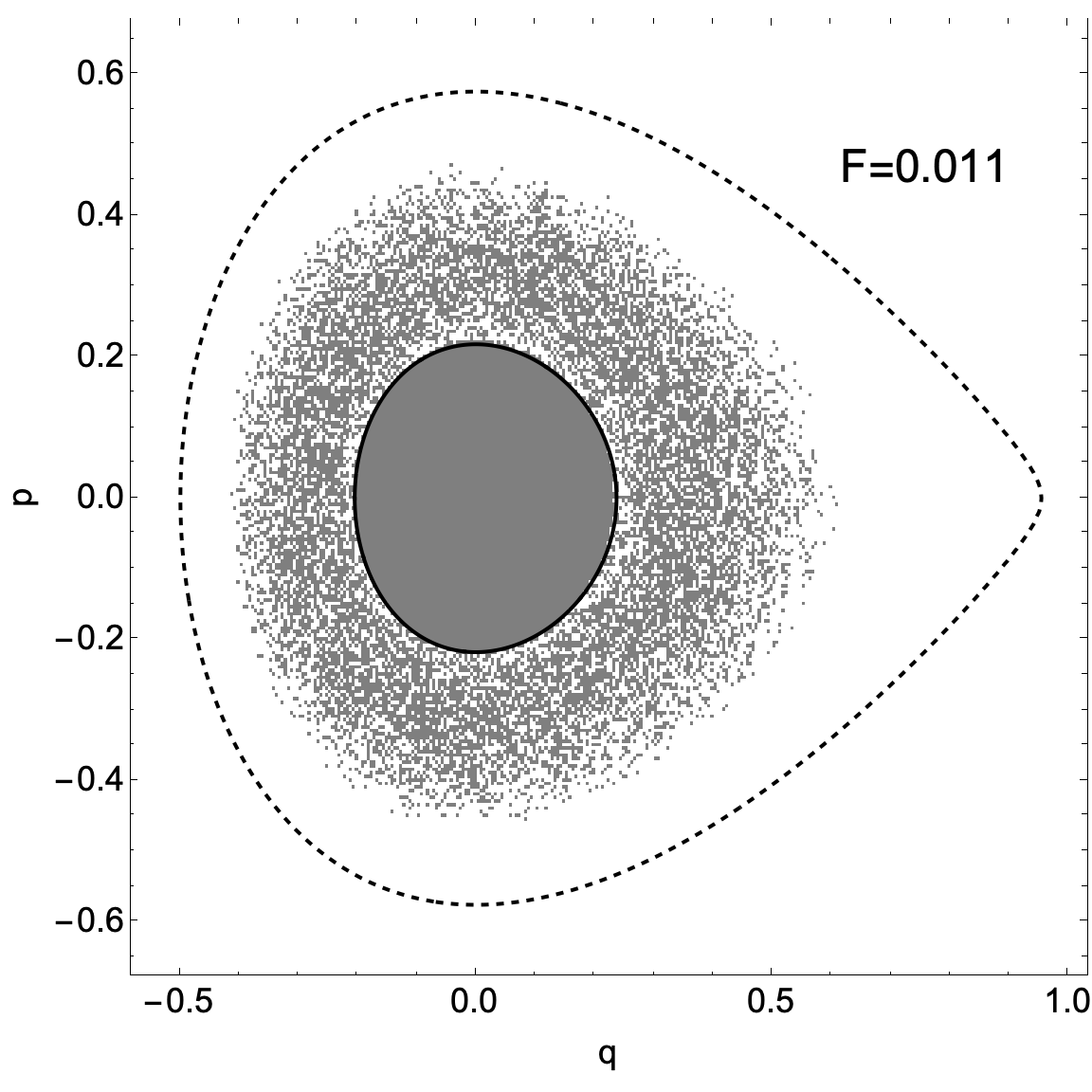}}
\hspace{10pt}
\subfloat[]{\includegraphics*[width=0.3\textwidth]{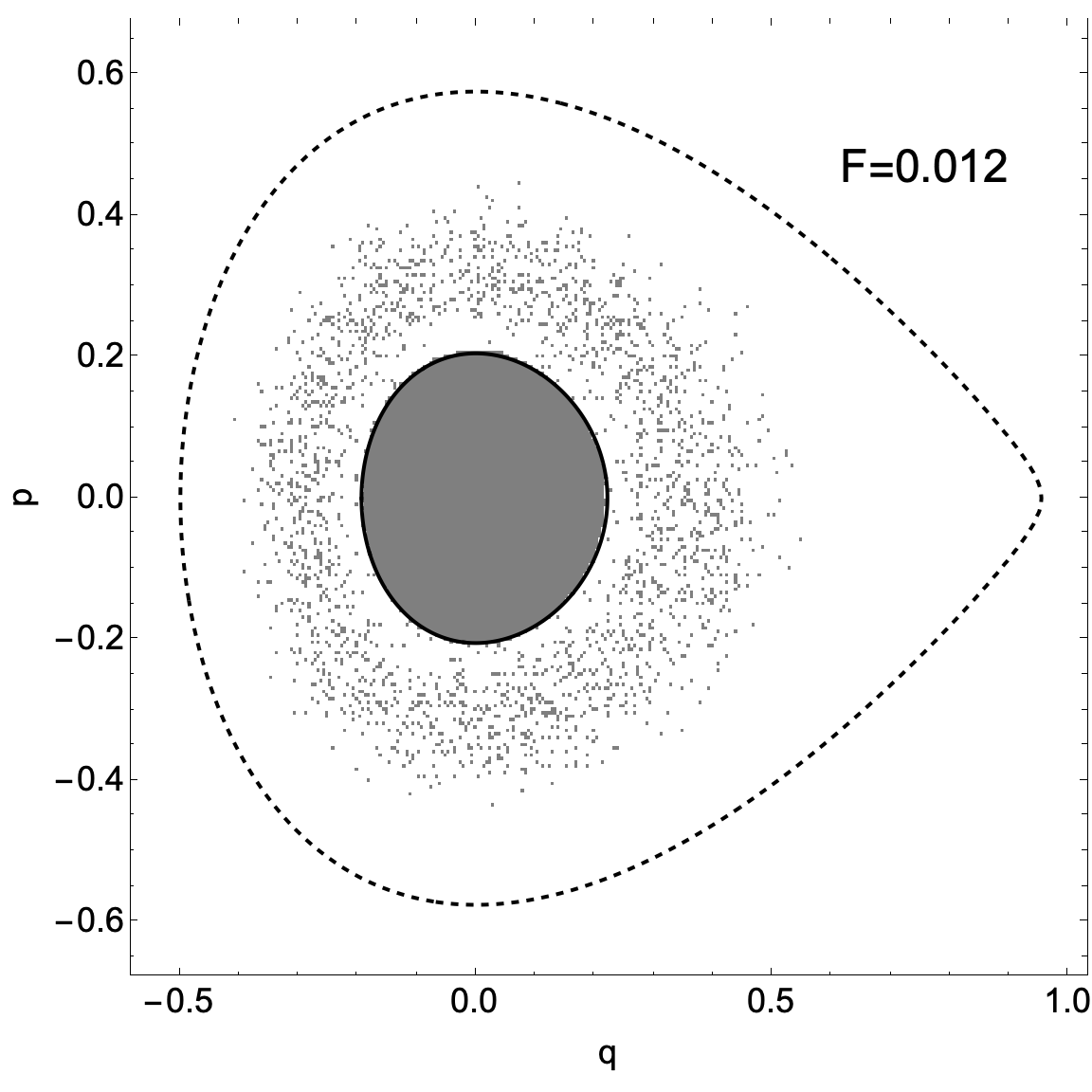}}
\hspace{10pt}
\subfloat[]{\includegraphics*[width=0.3\textwidth]{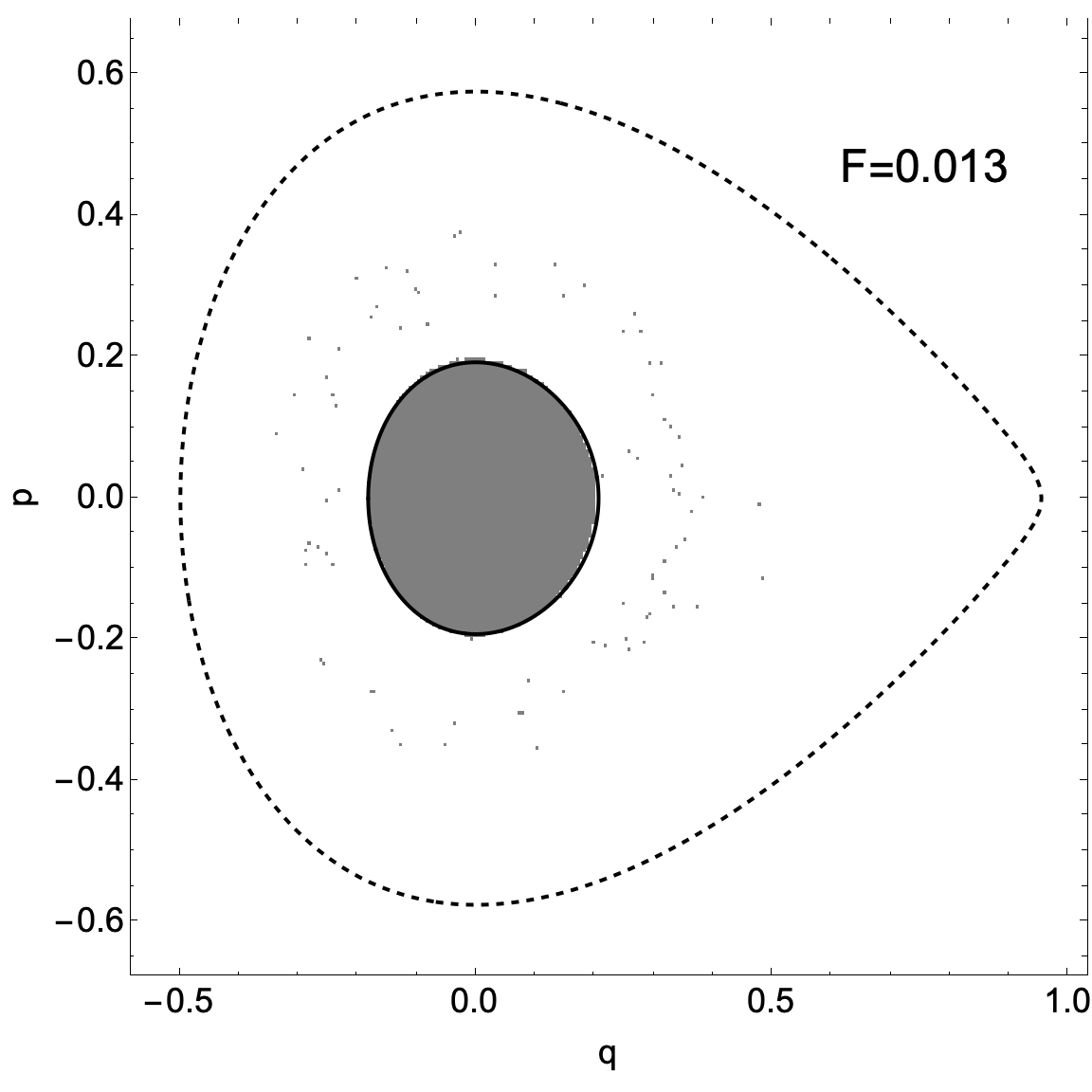}}\\
\subfloat[]{\includegraphics*[width=0.3\textwidth]{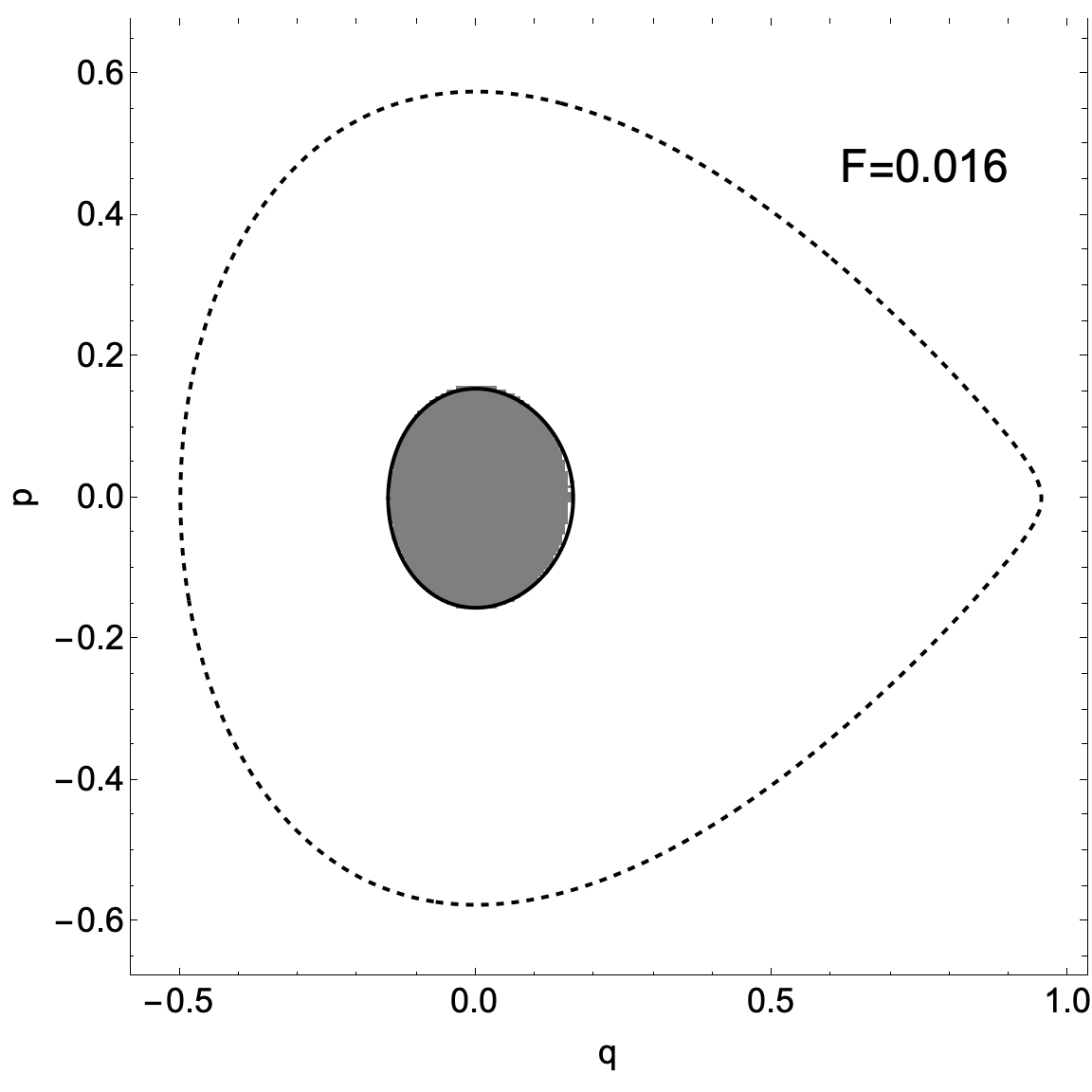}}
\hspace{10pt}
\subfloat[]{\includegraphics*[width=0.3\textwidth]{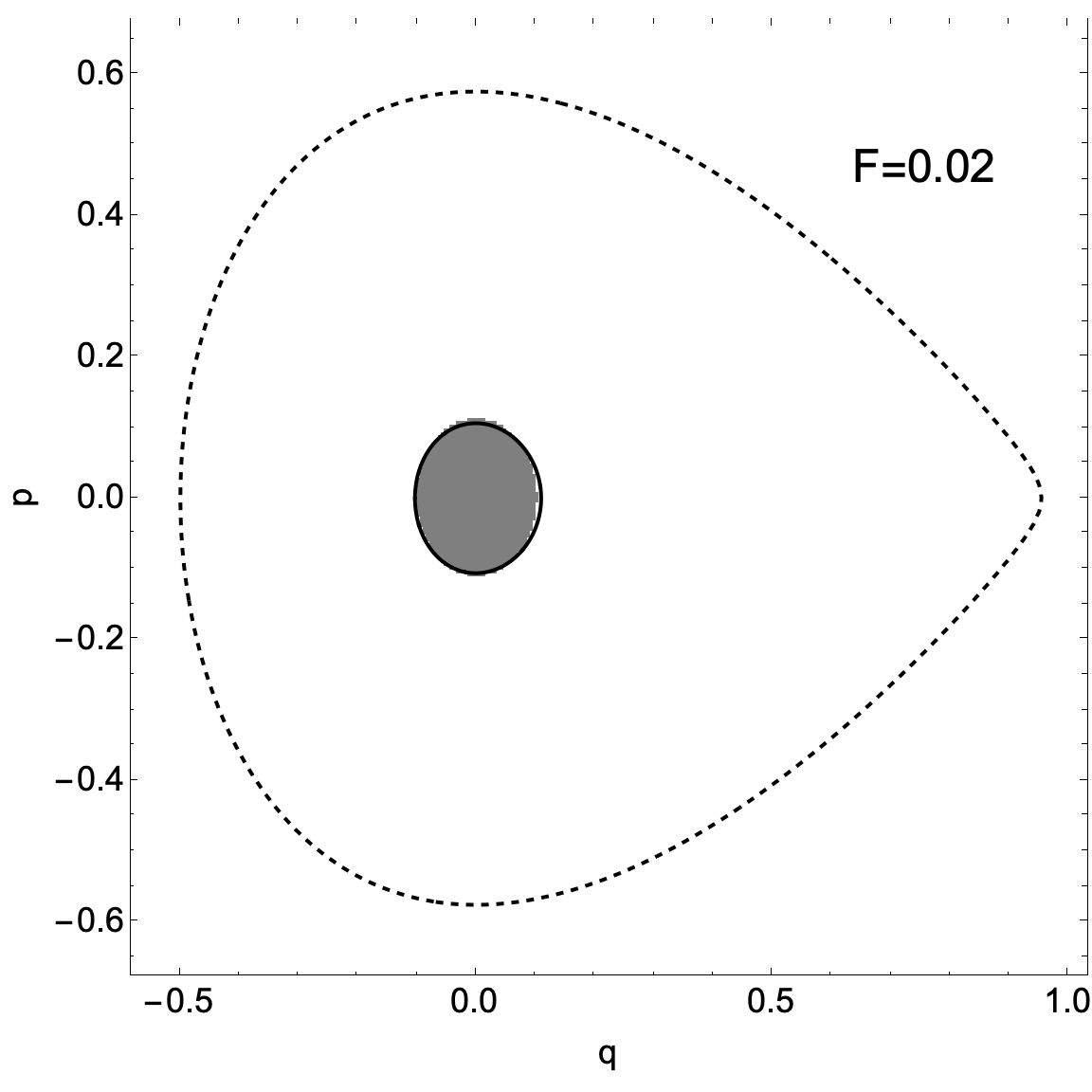}}
\hspace{10pt}
\subfloat[]{\includegraphics*[width=0.3\textwidth]{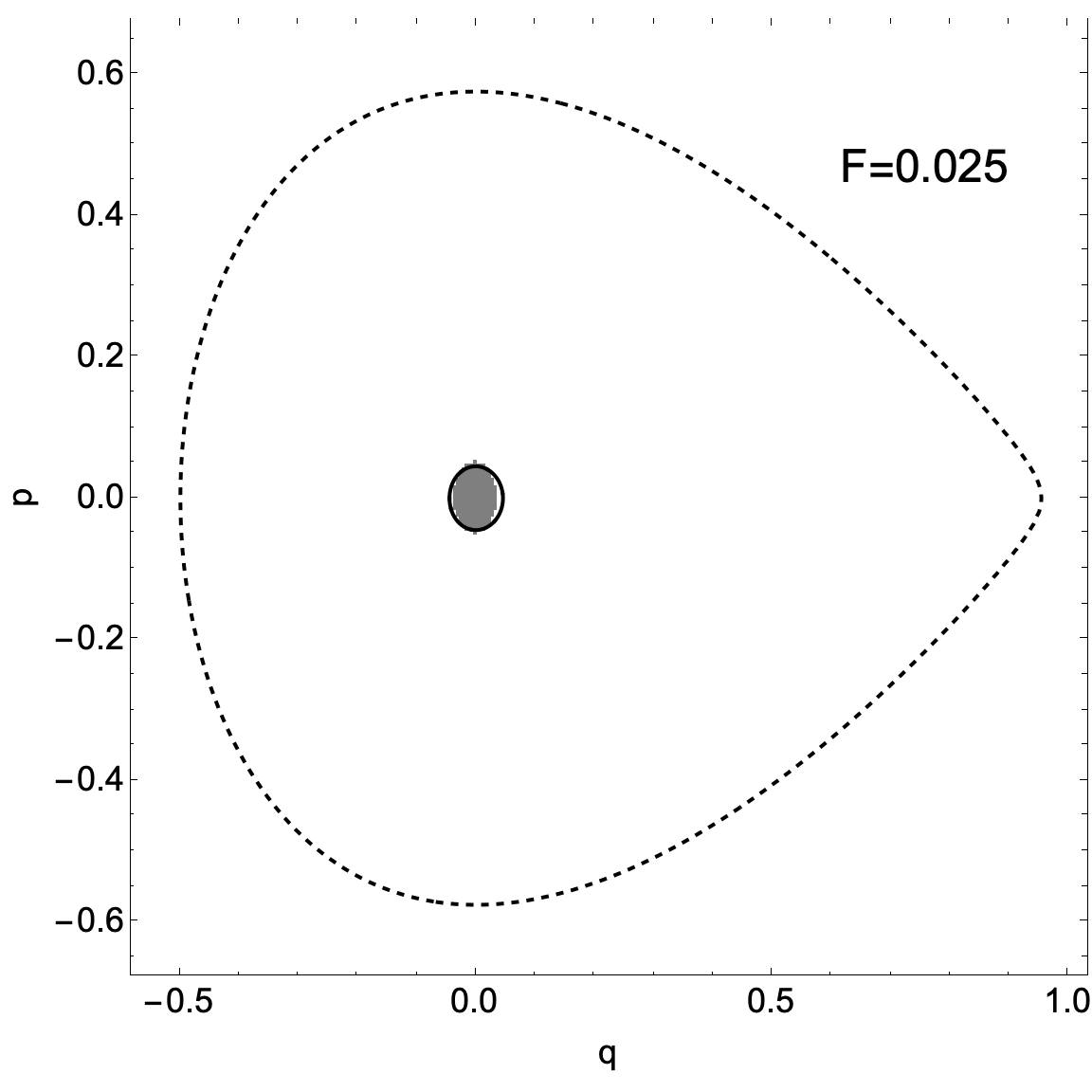}}
\caption{Evolution of the ``true" safe basin. Gray regions correspond to the numerical ``true" SB taken as an intersection of SBs for $\psi=0,\;\pi/11\dots 2\pi$. Black curve encircles analytic approximation of the ``true" SB. The value of frequency~$\Omega$ is the same as in Figure~\ref{fig:sb:evolution:phase:0}\label{fig:sb:evolution}}
\end{figure}

Figure~\ref{fig:sb:evolution} shows the evolution of the ``true" safe basins with increasing forcing amplitude~$F$. The numerical ``true" SB is obtained by simulating the SBs for each $\psi = 0,\; \pi/11,\; 2\pi/11,\dots, 2\pi$, and taking their intersection. As one can notice, the correspondence between the approximate SB and the numerics is remarkable for the most values of~$F$ except for the vicinity of $\widehat{F}$, i.e., when the SB collapses.

Figure~\ref{fig:true:sb:area:comparison} shows the comparison of the numerical and approximate erosion profiles of the ``true" SB. The largest discrepancy between the numerics and the approximation is observed around the sudden jump, i.e.,~$F\approx 0.01$. We attribute it to the fact that only~$21$ steps of the phase~$\psi$ were used. One can look at it from the following perspective: for the values of~$F$ close to but greater than~$\widehat{F}$, the SBMT is of kind I, i.e., the range of~$\vartheta$ does not span the whole period~$\left[0,\; 2\pi\right)$, however the ``gap" in~$\vartheta$ is narrow. Therefore, in order to obtain the numeric ``true" SB one has to take an intersection of SBs for significantly greater number of values~$\psi$. Figure~\ref{fig:gim:true:comparison:zoom} shows three points of the erosion profile computed using higher number of steps along the the phase~$\psi$ for the ``true" SB increases the accuracy, however, the convergence is notably slow.

\section{Conclusions and Discussions}\label{sec:conclusion}

In this research, we have provided significant insights into the study of safe basins and their erosion through the use of the Approximation of Isolated Resonance (AIR) method. By adopting the averaged energy and the slow phase as study variables, we have substantiated their effectiveness in the evaluation of safe basins.

Our work underscores the utility of the AIR method in determining the location and shape of the safe basins, as well as predicting their transformations with increasing value of the parameter~$F$. A key contribution is the successful application of the AIR method for the efficient exploration of the ``true" safe basins. Compared to numerical evaluations, which can be resource-intensive, our method offers a far more efficient approach.

However, we recognize that the current methodology is not without its limitations. A fundamental drawback is the method's reliance on analytic derivations of action-angle variables, which may not be straightforward or feasible for some complex systems. Moreover, the method fails to capture the stochastic layer near the separatrices, necessitating the use of numerical approximations for the appropriate truncation level. These limitations highlight areas for potential improvement in future studies.

In light of the identified drawbacks, an avenue for further research could be the estimation of the width of the stochastic layer. It is plausible that an accurate estimation technique might contribute to a more refined and comprehensive understanding of safe basins and their erosion.

In conclusion, while we recognize that improvements are needed, this research presents a pivotal step towards an efficient and effective method for studying safe basins. The insights offered by this study can serve as a robust foundation for further investigations aimed at refining the methodology and expanding its applications.

\section*{Acknowledgements}

Authors express their gratitude to the German Research Foundation (DFG) for the financial support within the project number 508244284. Authors also extend their heartfelt thanks to Attila Genda for his invaluable contributions and help in writing this manuscript.

\bibliographystyle{unsrt}
\bibliography{biblio.bib}

\begin{thebibliography}{10}

\bibitem{Virgin1989}
Lawrence~N. Virgin.
\newblock {Approximate criterion for capsize based on deterministic dynamics}.
\newblock {\em Dynamics and Stability of Systems}, 4(1), 1989.

\bibitem{belenky2007}
Vadim Belenky and Nikita Sevastianov.
\newblock {\em Stability and safety of ships: risk of capsizing}, pages
  165--289.
\newblock Society of Naval Architects and Marine Engineers, 2 edition, 2007.

\bibitem{mann2009energy}
BP~Mann.
\newblock Energy criterion for potential well escapes in a bistable magnetic
  pendulum.
\newblock {\em J. Sound Vib.}, 323(3-5):864--876, 2009.

\bibitem{5482087}
Fadi~M. Alsaleem, Mohammad~I. Younis, and Laura Ruzziconi.
\newblock An experimental and theoretical investigation of dynamic pull-in in
  mems resonators actuated electrostatically.
\newblock {\em Journal of Microelectromechanical Systems}, 19(4):794--806,
  2010.

\bibitem{leus2008dynamic}
Vitaly Leus and David Elata.
\newblock On the dynamic response of electrostatic mems switches.
\newblock {\em Journal of Microelectromechanical Systems}, 17(1):236--243,
  2008.

\bibitem{kramers1940}
H.A. Kramers.
\newblock Brownian motion in a field of force and the diffusion model of
  chemical reactions.
\newblock {\em Physica}, 7(4):284 -- 304, 1940.

\bibitem{fleming1993activated}
Graham~R Fleming and Peter Hanggi.
\newblock {\em Activated Barrier Crossing: applications in physics, chemistry
  and biology}, volume~4.
\newblock World Scientific, 1993.

\bibitem{VIRGIN1992357}
Lawrence~N. Virgin, Raymond~H. Plaut, and Ching-Chuan Cheng.
\newblock Prediction of escape from a potential well under harmonic excitation.
\newblock {\em International Journal of Non-Linear Mechanics}, 27(3):357--365,
  1992.

\bibitem{champneys2019happy}
Alan~R Champneys, Timothy~J Dodwell, Rainer~MJ Groh, Giles~W Hunt, Robin~M
  Neville, Alberto Pirrera, Amir~H Sakhaei, Mark Schenk, and M~Ahmer Wadee.
\newblock Happy catastrophe: recent progress in analysis and exploitation of
  elastic instability.
\newblock {\em Frontiers in Applied Mathematics and Statistics}, 5:34, 2019.

\bibitem{barone1982physics}
Antonio Barone and Gianfranco Paterno.
\newblock {\em Physics and applications of the Josephson effect}, pages
  136--160.
\newblock Wiley, 1982.

\bibitem{Gendelman2018}
Oleg Gendelman.
\newblock {Escape of a harmonically forced particle from an infinite-range
  potential well: a transient resonance}.
\newblock {\em Nonlinear Dynamics}, 93(1), 2018.

\bibitem{Karmi2021}
Gleb Karmi, Pavel Kravetc, and Oleg Gendelman.
\newblock {Analytic exploration of safe basins in a benchmark problem of forced
  escape}.
\newblock {\em Nonlinear Dynamics}, 106(3):1573--1589, 2021.

\bibitem{attila}
Attila Genda, Alexander Fidlin, and Oleg Gendelman.
\newblock {The level-crossing problem of a weakly damped particle in quadratic
  potential well under harmonic excitation}.
\newblock {\em Nonlinear Dynamics}, 111(22):20563--20578, 2023.

\bibitem{10.1063/5.0142761}
Pavel Kravetc, Oleg Gendelman, and Alexander Fidlin.
\newblock {Resonant escape induced by a finite time harmonic excitation}.
\newblock {\em Chaos: An Interdisciplinary Journal of Nonlinear Science},
  33(6), 06 2023.
\newblock 063116.

\bibitem{doi:10.1098/rspa.1989.0009}
John Michael~Tutill Thompson.
\newblock Chaotic phenomena triggering the escape from a potential well.
\newblock {\em Proceedings of the Royal Society of London. A. Mathematical and
  Physical Sciences}, 421(1861):195--225, 1989.

\bibitem{soliman1989integrity}
MS~Soliman and JMT Thompson.
\newblock Integrity measures quantifying the erosion of smooth and fractal
  basins of attraction.
\newblock {\em Journal of Sound and Vibration}, 135(3):453--475, 1989.

\bibitem{REGA2005902}
Giuseppe Rega and Stefano Lenci.
\newblock Identifying, evaluating, and controlling dynamical integrity measures
  in non-linear mechanical oscillators.
\newblock {\em Nonlinear Analysis: Theory, Methods \& Applications},
  63(5):902--914, 2005.
\newblock Invited Talks from the Fourth World Congress of Nonlinear Analysts
  (WCNA 2004).

\bibitem{Habib2021}
Giuseppe Habib.
\newblock {Dynamical integrity assessment of stable equilibria: a new rapid
  iterative procedure}.
\newblock {\em Nonlinear Dynamics}, 106(3):2073--2096, 2021.

\bibitem{lenci2003optimal}
Stefano Lenci and Giuseppe Rega.
\newblock Optimal control of homoclinic bifurcation: theoretical treatment and
  practical reduction of safe basin erosion in the helmholtz oscillator.
\newblock {\em Journal of Vibration and Control}, 9(3-4):281--315, 2003.

\bibitem{lenci2019global}
Stefano Lenci, Giuseppe Rega, et~al.
\newblock {\em Global Nonlinear Dynamics for Engineering Design and System
  Safety}, volume 588.
\newblock Springer, 2019.

\bibitem{rega2021global}
Giuseppe Rega and Valeria Settimi.
\newblock Global dynamics perspective on macro-to nano-mechanics.
\newblock {\em Nonlinear Dynamics}, 103(2):1259--1303, 2021.

\bibitem{regaLenci2008}
Giuseppe Rega and Stefano Lenci.
\newblock Dynamical integrity and control of nonlinear mechanical oscillators.
\newblock {\em Journal of Vibration and Control}, 14(1-2):159--179, 2008.

\bibitem{FARID2020105182}
Maor Farid.
\newblock Escape of a harmonically forced classical particle from asymmetric
  potential well.
\newblock {\em Communications in Nonlinear Science and Numerical Simulation},
  84:105182, 2020.

\end{thebibliography}

\end{document}